\documentclass[12pt]{amsart} 
\usepackage[active]{srcltx}
\usepackage{color}
\usepackage[english]{babel}
\usepackage[latin1]{inputenc}
\usepackage{amssymb}
\usepackage{amsmath}
\usepackage{amscd}
\usepackage{latexsym}
\usepackage{amsthm}
\usepackage{mathrsfs}

\theoremstyle{plain}
\newtheorem{thm}{Theorem}[section]
\newtheorem{cor}[thm]{Corollary}
\newtheorem{lem}[thm]{Lemma}
\newtheorem{prop}[thm]{Proposition}

\newtheorem{defn}[thm]{Definition}
\newtheorem{oss}[thm]{Remark}

\newcommand{\C}{\mathbb{C}} 

\renewcommand{\P}{\mathbb{P}}

\newcommand{\Z}{\mathbb{Z}}

\def\ch{\mbox{ch}}

\def\d{\mbox{d}} 
 
\def\Im{\mbox{Im}}

\def\Hom{\mbox{Hom}}

\def\td{\mbox{td}}

\begin{document}

\title{Moduli spaces of bundles over non-projective K3 surfaces}
\author[Perego, Toma]{Arvid Perego, Matei Toma}

\begin{abstract}
We study moduli spaces of sheaves over non-projective K3 surfaces. More precisely, if $v=(r,\xi,a)$ is a Mukai vector on a K3 surface $S$ with $r$ prime to $\xi$ and $\omega$ is a "generic" K\"ahler class on  $S$, we show that the moduli space $M$ of $\mu_{\omega}-$stable sheaves on $S$ with associated Mukai vector $v$ is an irreducible holomorphic symplectic manifold which is deformation equivalent to a Hilbert scheme of points on a K3 surface. If $M$ parametrizes only locally free sheaves, it is moreover hyperk\"ahler. Finally, we show that there is an isometry between $v^{\perp}$ and $H^{2}(M,\mathbb{Z})$ and that $M$ is projective if and only if $S$ is projective.
\end{abstract}

\date{\today}
\thanks{
 }

\maketitle
\tableofcontents

\section{Introduction}

Moduli spaces of sheaves on projective K3 surfaces have been studied since the '80s. In \cite{Fu} Fujiki considered the Hilbert scheme $Hilb^{2}(S)$ of 2 points on a K3 surface $S$; his result was widely generalized by Beauville in \cite{B}, who studied $Hilb^{n}(S)$ for any $n\in\mathbb{N}$, showing that it is an irreducible 
hyperk\"ahler manifold, i.e. a compact K\"ahler manifold which is simply connected, holomorphically symplectic and has $h^{2,0}=1$.

Moduli spaces of $\mu-$stable sheaves are a generalization of Hilbert schemes of points, and they have been extensively studied when the base surface $S$ is a projective K3 surface. In \cite{M1} Mukai showed that on the moduli space $M$ of simple sheaves of Mukai vector $v=(r,c_{1}(L),a)$ (i.e. of rank $r$, determinant $L$ and second Chern character $a-r$), there is a natural holomorphic symplectic form associated to the one on $S$. This moduli space $M$ is a non-separated scheme containing as a smooth open subset the moduli space $M^{\mu}_{v}(S,H)$ of $\mu_{H}-$stable sheaves (with respect to some ample line bundle $H$ on $S$) of Mukai vector $v$; Mukai's construction thus produces a holomorphic symplectic form on $M^{\mu}_{v}(S,\omega)$.

If $H$ is generic and $r$ and $L$ are prime to each other, then $M_{v}^{\mu}(S,H)$ is a projective holomorphically symplectic manifold. Moreover, it is an irreducible 
hyperk\"ahler manifold deformation equivalent to a Hilbert scheme of points on $S$ (see \cite{OG} and \cite{Y1}). 

If $S$ is a non-projective K3 surface and $\omega$ is a K\"ahler class on it, one still defines the notion of $\mu_{\omega}-$stable sheaf and constructs the moduli space $M^{\mu}_{v}(S,\omega)$ of $\mu_{\omega}-$stable sheaves of Mukai vector $v$. In \cite{T} it is shown that $M^{\mu}_{v}(S,\omega)$ is a smooth complex manifold carrying a holomorphic symplectic form. If $\omega$ is generic and $r$ is prime with $c_{1}(L)$, then $M^{\mu}_{v}(S,\omega)$ is even compact (see subsection \ref{genericity} for the precise notion of genericity we use for K\"ahler classes, 
called $v-$genericity in analogy to the projective case).

It is natural to ask if $M^{\mu}_{v}(S,\omega)$ is irreducible symplectic, and in this case what is its deformation class. We first show the following:

\begin{thm}
\label{thm:main}Let $S$ be a K3 surface, $v=(r,\xi,a)\in H^{*}(S,\mathbb{Z})$ where $\xi\in NS(S)$, $r>1$ prime with $\xi$ and $v^{2}\geq 0$. Suppose $\omega$ to be $v-$generic.
\begin{enumerate}
 \item The moduli space $M^{\mu}_{v}(S,\omega)$ is a compact, connected complex manifold of dimension $v^{2}+2$ which is holomorphically symplectic and deformation equivalent to a Hilbert scheme of points on a projective K3 surface.
 \item On $H^{2}(M^{\mu}_{v}(S,\omega),\mathbb{Z})$ there is a non-degenerate quadratic form, and there is an isometry between $H^{2}(M^{\mu}_{v},\mathbb{Z})$ and $v^{\perp}$ if $v^{2}>0$ (resp. $v^{\perp}/\mathbb{Z}v$ if $v^{2}=0$). 
\end{enumerate}
\end{thm}

The condition $v^{2}\geq 0$ implies that $M^{\mu}_{v}(S,\omega)\neq\emptyset$ (see \cite{BaLeP87}, \cite{TeTo}, \cite{KuYo}). As recalled above, if $S$ is projective and $\omega=c_{1}(H)$ for a generic ample line bundle $H$ we even know that $M^{\mu}_{v}(S,\omega)$ is an irreducible symplectic manifold. To prove Theorem \ref{thm:main}, we study the two remaining cases: $S$ is projective and $\omega\notin NS(S)$; and $S$ is non-projective.

When $S$ is projective and $\omega$ is not the first Chern class of an ample line bundle, we show that there is a $v-$generic ample line bundle $H$ such that $M^{\mu}_{v}(S,\omega)=M_{v}^{\mu}(S,H)$. This is done by showing that the $v-$chamber in which $\omega$ lies 
intersects the ample cone, and that moving the polarization inside a $v-$chamber does not affect the moduli space.

When $S$ is non-projective, the strategy to prove Theorem \ref{thm:main} is to deform $M^{\mu}_{v}(S,\omega)$ along the twistor family $\mathcal{X}\longrightarrow\mathbb{P}^{1}$ of $(S,\omega)$: even if the sheaves in $M^{\mu}_{v}(S,\omega)$ do not necessarily deform along such a twistor family, we can still deform them as twisted sheaves.

We then provide a construction of a relative moduli space of stable twisted sheaves extending Yoshioka's construction in \cite{Y3} to non-projective base manifolds and we show that we can connect the K3 surface $S$ to a projective K3 surface $S'$ only by means of twistor lines, in such a way that $M^{\mu}_{v}(S,\omega)$ deforms to $M^{\mu}_{v}(S',\omega')$ for some $v-$generic polarization $\omega'$ on $S'$. Theorem \ref{thm:main} holds true even if we replace $M^{\mu}_{v}(S,\omega)$ with a moduli space of stable twisted sheaves.
 
The non-degenerate quadratic form on $H^{2}(M^{\mu}_{v}(S,\omega),\mathbb{Z})$ is defined 
as a quadratic form on the second complex cohomology using the same definition of the Beauville form, the only difference being that we have to fix one holomorphic symplectic form to define it as a priori we have $h^{2,0}\geq 1$. We then show that it is non-degenerate. The construction of the isometry with $v^{\perp}$ is standard, and uses the same strategy as in the projective case.

As one might see from the statement on Theorem \ref{thm:main}, there is only one missing property for $M^{\mu}_{v}(S,\omega)$ to be an irreducible symplectic manifold; namely, we don't know if $M^{\mu}_{v}(S,\omega)$ is K\"ahler. This is a longstanding problem: on the open subset $M^{\mu-lf}_{v}(S,\omega)$ of $M^{\mu}_{v}(S,\omega)$ parametrizing locally free sheaves we have a natural K\"ahler metric, the Weil-Petersson metric, cf. \cite{Ito}, \cite{ItNa}, but at present nothing is known as to how this metric could extend to a K\"ahler metric on the whole $M^{\mu}_{v}(S,\omega)$.

The strategy to prove Theorem \ref{thm:main} together with \cite[Theorem 3.3]{HKLR} may be employed to obtain another proof of the existence of a K\"ahler metric on $M^{\mu-lf}_{v}(S,\omega)$ and of a description of a twistor family for such a hyperk\"ahler metric. But, as pointed to us by Daniel Huybrechts, this strategy does not allow to show that $M_{v}^{\mu}(S,\omega)$ carries a K\"ahler metric too. Let us remark however that there are choices of Mukai vectors for which $M_{v}^{\mu}(S,\omega)$ coincides with $M^{\mu-lf}_{v}(S,\omega)$ and is therefore a compact irreducible 
hyperk\"ahler manifold. Moreover such compact moduli spaces of stable locally free sheaves may acquire any positive even complex dimension; see Proposition \ref{Prop:dim}. 

As an application of the previous result, we will show the following projectivity criterion for the moduli spaces of slope-stable sheaves on a K3 surface:

\begin{thm}
\label{thm:proj}Let $S$ be a K3 surface, $v=(r,\xi,a)\in H^{2*}(S,\mathbb{Z})$ where $\xi\in NS(S)$, $r\geq 2$, $(r,\xi)=1$ and $v^{2}\geq 0$. If $\omega$ is a $v-$generic polarization, the moduli space $M^{\mu}_{v}(S,\omega)$ is projective if and only if $S$ is projective.
\end{thm}

\subsection*{Acknowledgements}
We are grateful to Daniel Huybrechts for pointing out to us a mistake in a previous version of this paper. We also thank the referee for his remarks and suggestions which significantly contributed to improve the exposition.

\section{Moduli spaces of stable sheaves}

In the following $S$ will be a K3 surface, possibly non-projective. If $\mathscr{F}$ is a coherent sheaf on $S$, we let the \textit{Mukai vector} of $\mathscr{F}$ be $$v(\mathscr{F}):=ch(\mathscr{F})\cdot\sqrt{td(S)}\in H^{2*}(S,\mathbb{Z}).$$If $v_{i}$ is the component of $v(\mathscr{F})$ in $H^{2i}(S,\mathbb{Z})$, we have $v_{0}=rk(\mathscr{F})$, $v_{1}=c_{1}(\mathscr{F})$ and $v_{2}=ch_{2}(\mathscr{F})+rk(\mathscr{F})=\frac{1}{2}c_1^2(\mathscr{F})-c_2(\mathscr{F})+rk(\mathscr{F})$, which will be viewed as an integer (i.e. we fix an isomorphism $H^{4}(S,\mathbb{Z})\simeq\mathbb{Z}$).
 
We recall that on $H^{2*}(S,\mathbb{Z})$ we have a pure weight-two Hodge structure and a lattice structure with respect to the Mukai pairing (see Definitions 6.1.5 and 6.1.11 of \cite{HL}): the obtained lattice will be referred to as \textit{Mukai lattice}, and we will write $v^{2}$ for the square of $v\in H^{2*}(S,\mathbb{Z})$ with respect to the Mukai pairing. Explicitly, $v^2=v_1^2-2v_0v_2$.

When $v_0\neq 0$ we define the \textit{discriminant} of $v$, or respectively of $\mathscr{F}$ in case $v=v(\mathscr{F})$, as
$$\Delta(v):=\frac{1}{2v_0^2}v^2+1,$$ 
This coincides with the definition of \cite{BaLeP87} for instance, where $$\Delta(\mathscr{F})=\Delta(v(\mathscr{F}))=
\frac{1}{rk(\mathscr{F})}\bigg(c_{2}(\mathscr{F})-\frac{rk(\mathscr{F})-1}{2rk(\mathscr{F})}c_1^2(\mathscr{F})\bigg).$$

\subsection{The stability condition}\label{subsection:stability}

Let $g$ be a K\"ahler metric on $S$ and $\omega$ the associated K\"ahler class, that will be called a \textit{polarization} on $S$. If $\mathscr{F}\in Coh(S)$ has positive rank, the \textit{slope} of $\mathscr{F}$ with respect to $\omega$ is $$\mu_{\omega}(\mathscr{F}):=\frac{c_{1}(\mathscr{F})\cdot\omega}{rk(\mathscr{F})}.$$

\begin{defn}
A torsion-free coherent sheaf $\mathscr{F}$ is $\mu_{\omega}-${\rm stable} if for every coherent subsheaf $\mathscr{E}\subseteq\mathscr{F}$ such that $0< rk(\mathscr{E})< rk(\mathscr{F})$ we have $\mu_{\omega}(\mathscr{E})<\mu_{\omega}(\mathscr{F})$. 
If $\mu_{\omega}(\mathscr{E})\leq\mu_{\omega}(\mathscr{F})$ for all such subsheaves $\mathscr{E}$, then we say that $\mathscr{F}$ is $\mu_{\omega}-${\rm semistable}.
\end{defn}

The family of $\mu_{\omega}-$stable sheaves of Mukai vector $v$ admits a moduli space $M^{\mu}_{v}(S,\omega)$. If $S$ is projective and $\omega$ is the first Chern class of an ample line bundle $H$, then $M^{\mu}_{v}(S,\omega)$ is the moduli space $M^{\mu}_{v}(S,H)$ of $\mu_{H}-$stable sheaves on $S$ with Mukai vector $v$. We have the following proposition dealing also with the non-projective case (see \cite{T}).

\begin{prop}
\label{prop:ssympl}Let $S$ be a K3 surface, $v\in H^{2*}(S,\mathbb{Z})$ a Mukai vector and $\omega$ a polarization on $S$. The moduli space $M^{\mu}_{v}(S,\omega)$ is a smooth, holomorphically symplectic manifold (possibly non-compact) and, if it is not empty, its dimension is $v^{2}+2$.
\end{prop}

In the following we will restrict to the case of those $M^{\mu}_{v}(S,\omega)$ which are non-empty and compact. We introduce in the next section some hypothesis on $v$ and $\omega$ under which $M^{\mu}_{v}(S,\omega)$ is compact. We now present a condition which guarantees its non-emptyness, and even the existence of a stable vector bundle with respect to any polarization.

Recall that over any non-algebraic surface there exist non-filtrable holomorphic rank two vector bundles (see \cite{BaLeP87}, \cite{Tdiss} p.18). By definition they do not admit coherent subsheaves of rank one, hence they are stable with respect to any polarization. 

We now extend this type of result to arbitrary rank in the case of K\"ahler surfaces. Following \cite{BaLeP87} we say that a coherent sheaf on the surface $S$ is \textit{irreducible} if its only coherent subsheaf of lower rank is the zero sheaf. In particular, an irreducible sheaf is stable with respect to any polarization. We have the following result, about the existence of locally free irreducible vector bundles.

\begin{prop}
\label{existence}
Let $S$ be a K\"ahler non-algebraic compact complex surface, $r$ a positive integer and $\xi\in NS(S)$. Then there exists a bound $b:=b(r,\xi)\in \mathbb{Z}$ depending on $r$ and on $\xi$ such that for any integer $c\geq b$ there is on $S$ an irreducible locally free sheaf $\mathscr{F}$ of rank $r$, $c_{1}(\mathscr{F})=\xi$ and $c_{2}(\mathscr{F})=c$. 
\end{prop}

\proof If $r=2$, a statement of this type is proved in \cite{BaLeP87} and in \cite{Tdiss} without the K\"ahler assumption. The idea there was to look at the versal deformation space of a rather arbitrary coherent sheaf $\mathscr{F}$ and show that if $c_{2}\gg 0$ then $\mathscr{F}$ must contain irreducible objects. For $r>2$ we shall this time consider deformations of suitably chosen coherent sheaves and make essential use of the fact that $S$ is K\"ahler. In this way we shall reduce ourselves to the argument used by B\u anic\u a and Le Potier in the case when the algebraic dimension of $S$ is zero, \cite[Th\'eor\`eme~5.3]{BaLeP87}.

We proceed by induction on $r$. The statement is trivial for $r=1$ and already proven for $r=2$. Let then $r\ge 3$ and suppose that the statement is true for rank $r-1$. Take an irreducible locally free sheaf $\mathscr{E}$ on $S$ of rank $r-1$, $c_{1}(\mathscr{E})=\xi$ and $c_{2}(\mathscr{E})=c$. Consider an irreducible component $B$ of the versal deformation space of $\mathscr{F}_0:=\mathscr{O}_{S}\oplus\mathscr{E}$ and the corresponding family $\mathscr{F}$ of coherent sheaves over $S\times B$. 

We shall check that if $c\gg 0$, the relative Douady space $D_{(X\times B)/B}(\mathscr{F},k)$ of flat quotients of rank $k$ of $\mathscr{F}$ over $B$ does not cover $B$ for $1\leq k\leq r-1$. Let $b:D_{(X\times B)/B}(\mathscr{F},k)\to B$ be the natural morphism and $Q\subset B$ a relatively compact subdomain of $B$ containing the origin $0\in B$. Fujiki proved in \cite{Fuj84} that any irreducible component of $b^{-1}(Q)$ is proper over $Q$. By another result of Fujiki in \cite{Fuj79}, there are countably many such components.
 
The idea is  to show by a dimension count that very general neighbours of $\mathscr{F}_0$ are not in the image of $D_{(X\times B)/B}(\mathscr{F},k)$ for $2\leq k\leq r-2$. Remark that if $\mathscr{F}_b$ is such a neighbour sitting in a short exact sequence $$0\to F'\to \mathscr{F}_b\to F''\to 0$$
with $F''$ torsion-free, then $F'$ and $F''$ are irreducible of different ranks, hence $\Hom (F',F'')=0=\Hom(F'', F')$. This remark makes the arguments in the proof of \cite[Th\'eor\`eme~5.3]{BaLeP87} work by replacing the corresponding inequality in loc. cit. Lemme 5.12. Hence our statement.\endproof

\subsection{The $v-$genericity for K\"ahler forms}\label{genericity}

Let $S$ be a K3 surface and $\mathscr{K}_{S}$ its K\"ahler cone, which is an open and convex cone in $H^{1,1}(S)$. For $v=(r,\xi,a)$ with $r\geq 2$ and $\xi\in NS(S)$, we define a system of hyperplanes in 
$H^{1,1}(S)$, which is locally finite in $\mathscr{K}_{S}$ and has the property that for any $\omega\in\mathscr{K}_{S}$ not lying on such hyperplanes, a torsion free sheaf $\mathscr{F}$ on $S$ with $v(\mathscr{F})=v$ is $\mu_{\omega}$-stable if and only if it is $\mu_{\omega}$-semistable. Polarizations verifying this will be called $v-$generic.

\subsubsection{The notion of $v-$genericity}

To start, let $S$ be any compact K\"ahler surface and fix $r,c_{2}\in\mathbb{Z}$, $c_{1}\in NS(S)$ and suppose $r>0$. We let $\tau:=(r,c_{1},c_{2})$, and if $\mathscr{F}\in Coh(S)$ of rank $r$ and Chern classes $c_{1}$ and $c_{2}$, we call $\tau$ the \textit{topological type} of $\mathscr{F}$. If $S$ is a K3 surface and $\mathscr{F}\in Coh(S)$ has Mukai vector $v=(r,\xi,a)$, its topological type is $\tau_{v}=(r,\xi,\xi^{2}/2+r-a)$. 

Notice that the discriminant $\Delta(\mathscr{F})$ only depends on the topological type of $\mathscr{F}$, hence we can talk about the \textit{discriminant} $\Delta(\tau)$ of $\tau$: more precisely, if $\tau=(r,c_{1},c_{2})$ then $$\Delta(\tau)=\frac{1}{r}\bigg(c_{2}-\frac{r-1}{2r}c_1^2\bigg).$$

We set $$W_{\tau}:=\{D\in NS(S)\,|\,-\frac{r^{4}}{2}\Delta(\tau)\leq D^{2}<0\}$$and for every $\alpha\in H^{1,1}(S)$ we write 
$$\alpha^{\perp}:=\{\beta\in H^{1,1}(S)\,|\,\alpha\cdot\beta=0\}.$$When $\alpha\neq 0$, the set $\alpha^{\perp}$ is a hyperplane in $H^{1,1}(S)$. Using the same argument of Lemma 4.C.2 of \cite{HL}, one shows that if $\beta\in H^{1,1}(S)$, then there is a open neighbourhood $U$ of $\beta$ in $H^{1,1}(S)$ such that $U\cap D^{\perp}\neq\emptyset$ for at most a finite number of $D\in W_{\tau}$. If the surface $S$ is K3, we will use the notation $W_{v}$ for $W_{\tau_{v}}$.

\begin{defn}
For every $D\in W_{\tau}$, the hyperplane $D^{\perp}\cap\mathscr{K}_{S}$ will be called $\tau-${\em wall in the K\"ahler cone of} $S$. A connected component of $\mathscr{K}_{S}\setminus\bigcup_{D\in W_{\tau}}D^{\perp}$ is an open convex cone called $\tau-${\em chamber in the K\"ahler cone of} $S$. A K\"ahler class in a $\tau-$chamber of $\mathscr{K}_{S}$ is called $\tau-${\em generic polarization}.
\end{defn}

If $S$ is a K3 surface and $v$ is a Mukai vector, we will call $v-$\textit{wall in the K\"ahler cone} (resp. $v-$\textit{chamber in the K\"ahler cone}, $v-$\textit{generic polarization}) a $\tau_{v}-$wall in the K\"ahler cone (resp. a $\tau_{v}-$chamber in the K\"ahler cone, a $\tau_{v}-$generic polarization). 

Recall that the ample cone of $S$ is $Amp(S)=\mathscr{K}_{S}\cap NS_{\mathbb{R}}(S)$ (where $NS_{\mathbb{R}}(S)=NS(S)\otimes\mathbb{R}$): if $S$ is a projective K3 surface and $\mathcal{C}\subseteq\mathscr{K}_{S}$ is a $v-$chamber in the K\"ahler cone of $S$, then $\mathcal{C}\cap NS_{\mathbb{R}}(S)$ is a $v-$chamber in the ample cone of $S$ in the usual terminology: if $H$ is an ample line bundle on $S$, then $c_{1}(H)$ is a $v-$generic polarization if and only if $H$ is $v-$generic as in \cite{HL}.

\subsubsection{Compactness of $M^{\mu}_{v}(S,\omega)$ when $\omega$ is $v-$generic}

Using the same proof as in the projective case (see Theorem 4.C.3 of \cite{HL}), we show that $v-$generic polarizations enjoy the above stated property concerning the existence of properly semistable sheaves. 

\begin{lem}
\label{lem:nonstable}Let $\omega$ be a K\"ahler class on a compact K\"ahler surface $S$, and $\mathscr{F}$ a $\mu_{\omega}-$semistable sheaf of topological type $\tau=(r,\xi,c_{2})$. Suppose that there is $\mathscr{E}\subseteq\mathscr{F}$ of rank $0<s<r$, first Chern class $\zeta$ and such that $\mu_{\omega}(\mathscr{E})=\mu_{\omega}(\mathscr{F})$. Then $D:=r\zeta-s\xi$ is such that $$-\frac{r^4}{2}\Delta(\tau)\leq D^{2}\leq 0,$$and $D^{2}=0$ if and only if $D=0$.
\end{lem}

\proof We can suppose that $\mathscr{E}$ is saturated, so that $\mathscr{G}:=\mathscr{F}/\mathscr{E}$ is torsion free, $\mu_{\omega}-$semistable and of rank $r-s$. Notice that as $\mu_{\omega}(\mathscr{E})=\mu_{\omega}(\mathscr{F})$, we have $D\cdot\omega=0$. As $\omega$ is a K\"ahler class, from the Hodge Index Theorem we then have $D^{2}\leq 0$, and $D^{2}=0$ if and only if $D=0$. We then just need to show that $D^{2}\geq-\frac{r^4}{2}\Delta(\tau)$.

By definition of the discriminant, it follows that $$\Delta(\mathscr{F})-\frac{s}{r}\Delta(\mathscr{E})-\frac{r-s}{r}\Delta(\mathscr{G})=-\frac{D^{2}}{2s(r-s)r^2}.$$
Now, recall that the Bogomolov inequality is surely satisfied by $\mathscr{E}$ and $\mathscr{G}$, so that $\Delta(\mathscr{E}),\Delta(\mathscr{G})\geq 0$. But this implies that $$-D^{2}\leq 2s(r-s)r^2\Delta(\mathscr{F})=2s(r-s)r^2\Delta(\tau)\leq\frac{r^4}{2}\Delta(\tau),$$and we are done.\endproof

Using the main result of \cite{T} we then get the following:

\begin{prop}
\label{prop:compact}Let $S$ be a K3 surface, $r\geq 2$ an integer and $\xi\in NS(S)$ such that $(r,\xi)=1$. Let $a\in\mathbb{Z}$, $v:=(r,\xi,a)$ and $\omega$ a $v-$generic polarization. If $M^{\mu}_{v}(S,\omega)\neq\emptyset$, then it is a smooth, compact, holomorphically symplectic manifold.
\end{prop}

\proof The statement follows from the main result of \cite{T} if $S$ is non-algebraic. When $S$ is projective we shall show in section 3.1 that there exists some integer ample class $H$ in the same $v$-chamber as $\omega$. The (semi)stability with respect to $\omega$ or with respect to $H$ will then come down to the same thing and $M^{\mu}_{v}(S,\omega)$ will coincide with the Gieseker moduli space $M_{v}(S,H)$ of $H-$semistable sheaves, which is known to be smooth, projective and holomorphically symplectic (see \cite{HL}).
\endproof

\section{Projective K3 surfaces with non-ample polarizations}

In this section we prove that if $S$ is a projective K3 surface, $v=(r,\xi,a)$ is a Mukai vector with $(r,\xi)=1$ and $\omega$ is a $v-$generic polarization, then $M^{\mu}_{v}(S,\omega)$ is an irreducible holomorphically symplectic manifold, deformation equivalent to a Hilbert scheme of points on $S$. 

\subsection{Changing polarization in a chamber}

We first show that changing polarization inside a chamber does not affect the moduli space. The following adaptation of Lemma 4.C.5 from \cite{HL} to the case of K\"ahler polarizations works also on K\"ahler manifolds; see \cite{GT} Lemma 6.2. 

\begin{lem}
\label{lem:openness}
Let $\omega$, $\omega'$ be two K\"ahler classes on a compact K\"ahler manifold $X$ and $\mathscr{F}$  be a torsion free sheaf on $X$ which is $\mu_{\omega}-$stable but not $\mu_{\omega'}-$stable. Denote by $$[\omega,\omega']:=\{\omega_{t}:=t\omega'+(1-t)\omega\,|\,t\in[0,1]\}$$ the segment from $\omega$ to $\omega'$.Then there is a K\"ahler class $\omega_t\in [\omega, \omega']$ such that $\mathscr{F}$ is properly $\mu_{\omega_t}-$semistable.
\end{lem}

As a consequence of this, changing the polarization inside a chamber does not affect the moduli space. This is well-known for $v-$generic ample line bundles, and requires the same proof. We let $M^{\mu}_{\tau}(S,\omega)$ be the moduli space of $\mu_{\omega}-$stable sheaves whose topological type is $\tau$. If $S$ is a K3 surface, then $M^{\mu}_{\tau_{v}}(S,\omega)=M^{\mu}_{v}(S,\omega)$

\begin{prop}
\label{prop:inside}Let $S$ be a smooth projective surface and $\tau=(r,\xi,c_{2})$ such that $r\geq 2$ and $\xi\in NS(S)$. Let $\mathcal{C}$ be a $\tau-$chamber in the K\"ahler cone of $S$, and $\omega,\omega'\in\mathcal{C}$. Then $M^{\mu}_{\tau}(S,\omega)=M^{\mu}_{\tau}(S,\omega')$.
\end{prop}

\proof We show that if $\mathscr{F}$ is a $\mu_{\omega}-$stable sheaf of topological type $\tau$, then it is $\mu_{\omega'}-$stable as well. Indeed, suppose that $\mathscr{F}$ is not $\mu_{\omega'}-$stable. By Lemma \ref{lem:openness} this implies that there is $\omega_{t}\in[\omega,\omega']$ such that $\mathscr{F}$ is properly $\mu_{\omega_{t}}-$semistable. Hence there is $\mathscr{E}\subseteq\mathscr{F}$ of rank $0<s<r$ and first Chern class $\zeta$, such that $\mu_{\omega_t}(\mathscr{E})=\mu_{\omega_t}(\mathscr{F})$. 

Let $D:=r\zeta-s\xi$: hence $D\cdot\omega_{t}=0$, and by Lemma \ref{lem:nonstable} we have $D\in W_{\tau}\cup\{0\}$. Notice that as $\mathscr{F}$ is $\mu_{\omega}-$stable, we have $D\cdot\omega<0$, so that $D\in W_{\tau}$. It follows that $\omega_{t}\notin\mathcal{C}$ which is not possible as $\mathcal{C}$ is convex. In conclusion, $\mathscr{F}$ is $\mu_{\omega'}-$stable.\endproof

\subsection{Conclusion for projective K3 surfaces}

We first introduce some notations: if $S$ is a projective surface, we let $NS_{\mathbb{R}}(S)$ be the real Néron-Severi space of $S$, which is a linear subspace of $H^{1,1}(S)$. Recall that on $H^{1,1}(S)$ we have a non-degenerate intersection product whose restriction to $NS_{\mathbb{R}}(S)$ remains non-degenerate. Let $T_{\mathbb{R}}(S)$ be the orthogonal of $NS_{\mathbb{R}}(S)$ in $H^{1,1}(S)$, so that we have $H^{1,1}(S)=NS_{\mathbb{R}}(S)\oplus T_{\mathbb{R}}(S)$.

Finally, for every $\alpha\in H^{1,1}(S)$ we let $p_{NS}:H^{1,1}(S)\longrightarrow NS_{\mathbb{R}}(S)$ and $p_{T}:H^{1,1}(S)\longrightarrow T_{\mathbb{R}}(S)$ be the two projections. Moreover, for every $\alpha\in H^{1,1}(S)$ we let $\alpha_{NS}:=p_{NS}(\alpha)$ and $\alpha_{T}:=p_{T}(\alpha)$.

The first result we show is the following:

\begin{lem}
\label{lem:projomega}
Let $S$ be a projective surface and $\omega$ a K\"ahler class on $S$. 
\begin{enumerate}
 \item The class $\omega_{NS}$ is an ample class on $S$.
 \item For every $\xi\in NS_{\mathbb{R}}(S)$ we have $\xi\cdot\omega=\xi\cdot\omega_{NS}$.
\end{enumerate}
\end{lem}

\proof Recall that $\omega=\omega_{NS}+\omega_{T}$: it follows that
for every non-zero effective curve class $C$ we have
$$\omega_{NS}\cdot C=\omega\cdot C-\omega_{T}\cdot C=\omega\cdot C>0,$$since $\omega_{T}$ is orthogonal to $NS_{\mathbb{R}}(S)$ (where $C$ lies), and $\omega$ is a K\"ahler class. This implies that $\omega_{NS}$ is a nef class on $S$.

In particular, this means that $\omega_{NS}$ is a class in the closure of the ample cone of $S$. Now, recall that the projection $p_{NS}$ is an open map; moreover, the previous part of the proof shows that the image of the K\"ahler cone of $S$ under $p_{NS}$ is contained in the nef cone of $S$.

As the K\"ahler cone is open in $H^{1,1}(S)$ and the interior of the nef cone is the ample cone, it follows that the image of the K\"ahler cone by projection is contained in the ample cone.

The last point of the statement is simply the fact that $\omega_{T}$ is orthogonal to $NS_{\mathbb{R}}(S)$.\endproof

Using the previous Lemma, we can finally prove the following, which shows part (1) of Theorem \ref{thm:main}.

\begin{thm}
\label{thm:change}Let $S$ be a projective K3 surface and $v=(r,\xi,a)\in H^{2*}(S,\mathbb{Z})$ such that $r\geq 2$, $\xi\in NS(S)$ and $(r,\xi)=1$. If $\omega$ is $v-$generic and $M^{\mu}_{v}(S,\omega)\neq\emptyset$, then $M^{\mu}_{v}(S,\omega)$ is a projective irreducible hyperk\"ahler manifold deformation equivalent to a Hilbert scheme of points on $S$.
\end{thm}

\proof The class $\omega_{NS}$ is ample by Lemma \ref{lem:projomega}, and $\omega_{NS}\cdot\xi=\omega\cdot\xi$ for every $\xi\in NS_{\mathbb{R}}(S)$. It follows that for every $\mathscr{F}\in Coh(S)$ we have $\mu_{\omega}(\mathscr{F})=\mu_{\omega_{NS}}(\mathscr{F})$. In particular, a coherent sheaf is $\mu_{\omega}-$stable if and only if it is $\mu_{\omega_{NS}}-$stable, so that $M^{\mu}_{v}(S,\omega)=M^{\mu}_{v}(S,\omega_{NS})$.

Moreover, if $D\in W_{v}$, then $\omega_{NS}\cdot D=\omega\cdot D$: as $\omega$ is $v-$generic, it follows that $\omega_{NS}$ is $v-$generic. Let $\mathcal{C}$ be the $v-$chamber of the ample cone where $\omega_{NS}$ lies. As $\mathcal{C}$ is open in $Amp(S)$, there is $\epsilon>0$ such that the ball $B_{\epsilon}(\omega_{NS})\subseteq Amp(S)$ of ray $\epsilon$ and centred at $\omega_{NS}$ is contained in $\mathcal{C}$. Let $\omega'\in B_{\epsilon}(\omega_{NS})\cap H^{2}(S,\mathbb{Q})$: by Proposition \ref{prop:inside} we have $M^{\mu}_{v}(S,\omega_{NS})=M^{\mu}_{v}(S,\omega')$. 

As $\omega'\in H^{2}(S,\mathbb{Q})\cap H^{1,1}(S)$, there are $p\in\mathbb{N}$ and $H\in Pic(S)$ such that $p\omega'=c_{1}(H)$. As $\omega'\in\mathcal{C}$ and $\mathcal{C}$ is a cone, we have $c_{1}(H)\in\mathcal{C}$: hence $H$ is a $v-$generic ample line bundle, and $M^{\mu}_{v}(S,\omega')=M^{\mu}_{v}(S,H)$. By \cite{OG} and \cite{Y1} $M^{\mu}_{v}(S,H)$ is an irreducible hyperk\"ahler manifold deformation equivalent to a Hilbert scheme of points, and we are done.\endproof

\begin{oss}
\label{oss:intcham}
A useful corollary of Lemma \ref{lem:projomega} is that if $\mathcal{C}$ is a $v-$chamber in the K\"ahler cone of $S$, then $\mathcal{C}$ intersects the ample cone (and the intersection is clearly a $v-$chamber in the ample cone of $S$). 

Indeed, consider the segment $[\omega,\omega_{NS}]$: as the projection $p_{NS}$ is a linear map, we have that $[\omega,\omega_{NS}]\cap NS_{\mathbb{R}}(S)=\{\omega_{NS}\}$, and that $p_{NS}([\omega,\omega_{NS}])=\{\omega_{NS}\}$.
                                                                                                    
We show that $\omega_{NS}\in\mathcal{C}$ (showing that $\mathcal{C}$ intersects the ample cone by Lemma \ref{lem:projomega}). Indeed, suppose that $\omega_{NS}$ does not lie in $\mathcal{C}$: it follows that there is $\omega'\in[\omega,\omega_{NS}]$ lying on a $v-$wall. This means that there is $D\in W_{v}$ such that $\omega'\cdot D=0$. But as $p_{NS}(\omega')=p_{NS}(\omega)=\omega_{NS}$, it follows that $\omega\cdot D=0$, which is not possible.
\end{oss}

\section{Moduli spaces of stable twisted sheaves}

In this section we recall the notion of twisted sheaf on a complex manifold, and we introduce the notion of stability for coherent twisted sheaves. We will then construct (relative) moduli spaces of stable twisted sheaves on a K3 surface (not necessarily projective): they will be used to show that the moduli spaces $M^{\mu}_{v}(S,\omega)$ of $\mu_{\omega}-$stable sheaves with Mukai vector $v=(r,\xi,a)$ such that $r$ and $\xi$ are prime to each other are compact, connected, simply connected and deformation equivalent to a Hilbert scheme of points on a projective K3 surface (whenever the polarization $\omega$ is $v-$generic). 

\subsection{Twisted sheaves and stability}

We recall some basic facts about twisted sheaves on a complex manifold $X$ (we refer the interested reader to \cite{Cal} or \cite{L} for more details). 

There are several definitions of twisted sheaves, giving equivalent categories. We use three of them: the first one is due to C\u ald\u araru \cite{Cal}, and presents twisted sheaves as a twisted gluing of local coherent sheaves on $X$; the second one (to be found again in \cite{Cal}) presents twisted sheaves as modules over an Azumaya algebra on $X$; the last one, due to Yoshioka \cite{Y3}, presents twisted sheaves as a full subcategory of the category of coherent sheaves on some projective bundle over $X$. 

We begin by recalling these definitions. As our aim are moduli spaces of stable twisted sheaves, we need to introduce several notions: first, we recall the Chern character and the slope of a twisted sheaf (for projective manifolds, this was done in \cite{HS} and \cite{Y3}); we then introduce $\mu_{\omega}-$stability for twisted sheaves (with respect to a K\"ahler form $\omega$).  

\subsubsection{Twisted sheaves following C\u ald\u araru}

Let $\mathscr{U}=\{U_{i}\}_{i\in I}$ be an open covering of $X$, and let $U_{ij}:=U_{i}\cap U_{j}$ and $U_{ijk}:=U_{i}\cap U_{j}\cap U_{k}$. Choose a $2-$cocyle $\{\alpha_{ijk}\}$, where $\alpha_{ijk}\in\mathscr{O}_{X}^{*}(U_{ijk})$, defining a class $\alpha\in H^{2}(X,\mathscr{O}_{X}^{*})$. A $\{\alpha_{ijk}\}-$\textit{twisted coherent sheaf} is a collection $\mathscr{F}=\{\mathscr{F}_{i},\phi_{ij}\}$, where $\mathscr{F}_{i}\in Coh(U_{i})$ for every $i\in I$, and for every $i,j\in I$ $\phi_{ij}:\mathscr{F}_{j|U_{ij}}\longrightarrow\mathscr{F}_{i|U_{ij}}$ is an isomorphism in $Coh(U_{ij})$ such that 
\begin{enumerate}
 \item $\phi_{ii}=id_{\mathscr{F}_{i}}$ for every $i\in I$;
 \item $\phi_{ij}=\phi_{ji}^{-1}$ for every $i,j\in I$;
 \item $\phi_{ij}\circ\phi_{jk}\circ\phi_{ki}=\alpha_{ijk}\cdot id$ for every $i,j,k\in I$.
\end{enumerate}
By a \textit{morphism of $\{\alpha_{ijk}\}-$twisted sheaves} $$f:\mathscr{F}=\{\mathscr{F}_{i},\phi_{ij}\}\longrightarrow\mathscr{G}=\{\mathscr{G}_{i},\psi_{ij}\}$$we mean a collection $f=\{f_{i}\}$ of morphisms $f_{i}:\mathscr{F}_{i}\longrightarrow\mathscr{G}_{i}$ of $\mathscr{O}_{U_{i}}-$modules such that $\psi_{ij}\circ f_{j}=f_{i}\circ\phi_{ij}$ for every $i,j\in I$.

The $\{\alpha_{ijk}\}-$twisted coherent sheaves form an abelian category $Coh(X,\{\alpha_{ijk}\})$. If $\{\alpha_{ijk}\}$ and $\{\alpha'_{ijk}\}$ are two representatives of the same class $\alpha\in H^{2}(X,\mathscr{O}_{X}^{*})$, then there is an equivalence between $Coh(X,\{\alpha_{ijk}\})$ and $Coh(X,\{\alpha'_{ijk}\})$, so that we can speak of the category $Coh(X,\alpha)$ of coherent $\alpha-$twisted sheaves.

If $\mathscr{F}\in Coh(X,\alpha)$ and $\mathscr{G}\in Coh(X,\beta)$, we can define in a natural way $\mathscr{F}\otimes\mathscr{G}$ and $\mathscr{H}om(\mathscr{F},\mathscr{G})$: the first one is a coherent $\alpha\beta-$twisted sheaf, while the second is a coherent $\alpha^{-1}\beta-$twisted sheaf.

We now recall an important definition: a sheaf $\mathcal{A}$ of $\mathscr{O}_{X}-$modules is said to be an \textit{Azumaya algebra} if it is a sheaf of $\mathscr{O}_{X}-$algebras whose generic fibre is a central simple algebra. Equivalence classes of Azumaya algebras form a group $Br(X)$, the \textit{Brauer group} of $X$, which has an injection into $H^{2}(X,\mathscr{O}_{X}^{*})$. One of the main properties we will use in the following is (see Theorem 1.3.5 of \cite{Cal}):

\begin{prop}
\label{prop:exivbtw}
Let $X$ be a complex manifold and $\alpha\in Br(X)$. Then there exist a locally free $\alpha-$twisted sheaf on $X$.
\end{prop}

For the rest of this section, we suppose $\alpha\in Br(X)$ and define the twisted Chern character and twisted Mukai vector for $\alpha-$twisted sheaves. More precisely, let $\mathscr{F}$ be an $\alpha-$twisted coherent sheaf on $X$ and $E$ a locally free $\alpha-$twisted coherent sheaf. The \textit{Chern character} of $\mathscr{F}$ with respect to $E$ is $$ch_{E}(\mathscr{F}):=\frac{ch(\mathscr{F}\otimes E^{\vee})}{\sqrt{ch(E\otimes E^{\vee})}}.$$The \textit{Mukai vector} of $\mathscr{F}$ with respect to $E$ is $$v_{E}(\mathscr{F}):=ch_{E}(\mathscr{F})\cdot\sqrt{td(X)}.$$The \textit{slope} of a torsion-free $\alpha-$twisted sheaf $\mathscr{F}$ with respect to $E$ and to a K\"ahler class $\omega$ is $$\mu_{E,\omega}(\mathscr{F}):=\frac{c_{E,1}(\mathscr{F})\cdot\omega}{rk(\mathscr{F})},$$where $c_{E,1}(\mathscr{F})$ is the component of $ch_{E}(\mathscr{F})$ lying in $H^{2}(S,\mathbb{Q})$.

We collect some explicit formulas when $X=S$ is a K3 surface. Let $r:=rk(\mathscr{F})$, $s:=rk(E)$, $\xi:=c_{1}(\mathscr{F}\otimes E^{\vee})$, $a:=ch_{2}(\mathscr{F}\otimes E^{\vee})$ and $b:=ch_{2}(E\otimes E^{\vee})$. Then $$ch_{E}(\mathscr{F})=(r,\xi/s,(2as-rb)/2s^{2}),$$ $$v_{E}(\mathscr{F})=(r,\xi/s,r+(2as-rb)/2s^{2})$$so that
\begin{equation}
\label{eq:Eslope}
\mu_{E,\omega}(\mathscr{F})=\frac{\xi\cdot\omega}{rs}=\frac{c_{1}(\mathscr{F}\otimes E^{\vee})\cdot\omega}{rk(\mathscr{F})rk(E)}=\mu_{\omega}(\mathscr{F}\otimes E^{\vee})
\end{equation}
and 
\begin{equation}
\label{eq:vsquare}
v_{E}^{2}(\mathscr{F})=\frac{\xi^{2}}{s^{2}}-\frac{2ra}{s}+\frac{r^{2}b}{s^{2}}-2r^{2}.
\end{equation}
If $\alpha=0$, then one easily sees that $\mu_{E,\omega}(\mathscr{F})=\mu_{\omega}(\mathscr{F})-\mu_{\omega}(E)$ and that \begin{equation}
\label{eq:v2}
v^{2}_{E}(\mathscr{F})=v^{2}(\mathscr{F}).
\end{equation}
If $\mathscr{F}$ is a torsion free $\alpha-$twisted sheaf on $S$, we let $$ch_{\alpha}(\mathscr{F}):=ch_{\mathscr{F}^{\vee\vee}}(\mathscr{F}),\,\,\,\,\,\,\,\,\,\,v_{\alpha}(\mathscr{F}):=v_{\mathscr{F}^{\vee\vee}}(\mathscr{F}),$$called \textit{twisted Chern character} and \textit{twisted Mukai vector} of $\mathscr{F}$. The \textit{twisted slope} of $\mathscr{F}$ with respect to $\omega$ is $$\mu_{\alpha,\omega}(\mathscr{F}):=\frac{c_{\alpha,1}(\mathscr{F})\cdot\omega}{rk(\mathscr{F})},$$where $c_{\alpha,1}(\mathscr{F})$ is the component of $ch_{\alpha}(\mathscr{F})$ in $H^{2}(S,\mathbb{Q})$.

Using twisted slopes, we introduce the notion of stability for twisted sheaves. Fix $\alpha\in Br(X)$ and $E$ an $\alpha-$twisted locally free sheaf.

\begin{defn}
We say that a torsion-free $\mathscr{F}\in Coh(X,\alpha)$ is $\mu_{E,\omega}-${\em stable} if for every $\alpha-$twisted coherent subsheaf $\mathscr{E}\subseteq\mathscr{F}$ such that $0<rk(\mathscr{E})<rk(\mathscr{F})$ we have $\mu_{E,\omega}(\mathscr{E})<\mu_{E,\omega}(\mathscr{F})$. If $\mu_{E,\omega}(\mathscr{E})\leq\mu_{E,\omega}(\mathscr{F})$ for every such subsheaf, we say that $\mathscr{F}$ is $\mu_{E,\omega}-${\em semistable}. A $\mu_{\mathscr{F}^{\vee\vee},\omega}-$(semi)stable sheaf will be called $\mu_{\alpha,\omega}-${\em (semi)stable}.
\end{defn}

To conclude this section, we show that the $\mu_{E,\omega}-$stability does not depend on $E$.

\begin{lem}
\label{lem:properties}
Let $\alpha\in Br(S)$, $\mathscr{F}\in Coh(S,\alpha)$ and $\omega\in\mathscr{K}_{S}$. If $E',E\in Coh(S,\alpha)$ are locally free, then $\mathscr{F}$ is $\mu_{E,\omega}-$stable if and only if it is $\mu_{E',\omega}-$stable. In particular, $\mathscr{F}$ is $\mu_{E,\omega}-$stable if and only if it is $\mu_{\alpha,\omega}-$stable. If $\alpha=0$, the sheaf $\mathscr{F}$ is $\mu_{0,\omega}-$stable if and only if it is $\mu_{\omega}-$stable.
\end{lem}

\proof Let $\mathscr{F}\in Coh(S,\alpha)$, $\mathscr{G}$ an $\alpha-$twisted coherent subsheaf of $\mathscr{F}$, and $H$ a locally free $\alpha-$twisted coherent sheaf. Then $$rk(H)rk(\mathscr{F})c_{1}(\mathscr{G}\otimes H^{\vee})-rk(H)rk(\mathscr{G})c_{1}(\mathscr{F}\otimes H^{\vee})=$$ 
\begin{equation}
\label{eq:c1tens}
=c_{1}(\mathscr{G}\otimes\mathscr{F}^{\vee}\otimes H\otimes H^{\vee})=rk^{2}(H)c_{1}(\mathscr{G}\otimes\mathscr{F}^{\vee}).
\end{equation}

Suppose now that $\mathscr{F}$ is $\mu_{E,\omega}-$stable but not $\mu_{E',\omega}-$stable. Hence there is an $\alpha-$twisted coherent subsheaf $\mathscr{G}$ of $\mathscr{F}$ of rank $0<s<rk(\mathscr{F})$ such that $\mu_{E',\omega}(\mathscr{G})\geq\mu_{E',\omega}(\mathscr{F})$. By $\mu_{E,\omega}-$stability of $\mathscr{F}$ we even have $\mu_{E,\omega}(\mathscr{G})<\mu_{E,\omega}(\mathscr{F})$. Writing these two inequalities explicitly we have
\begin{equation}
\label{eq:disugprimo} 
\omega\cdot(rk(E')rc_{1}(\mathscr{G}\otimes(E')^{\vee})-rk(E')sc_{1}(\mathscr{F}\otimes(E')^{\vee}))\geq 0,
\end{equation}
\begin{equation}
\label{eq:disug}
\omega\cdot(rk(E)rc_{1}(\mathscr{G}\otimes E^{\vee})-rk(E)sc_{1}(\mathscr{F}\otimes E^{\vee}))< 0.
\end{equation}
Using equation (\ref{eq:c1tens}) for $H=E'$, equation (\ref{eq:disugprimo}) becomes $\omega\cdot c_{1}(\mathscr{G}\otimes\mathscr{F}^{\vee})\geq 0$. Using equation (\ref{eq:c1tens}) for $H=E$, equation (\ref{eq:disug}) becomes $\omega\cdot c_{1}(\mathscr{G}\otimes\mathscr{F}^{\vee})<0$, getting a contradiction.\endproof

\subsubsection{Twisted sheaves as $\mathcal{A}-$modules}

Let again $X$ be a complex manifold and $\mathcal{A}$ an Azumaya algebra on $X$. We let $Coh(X,\mathcal{A})$ be the abelian category of coherent sheaves on $X$ having the structure of $\mathcal{A}-$module. The following is Proposition 1.3.6 of \cite{Cal}:

\begin{prop}
\label{prop:equivazu}Let $X$ be a complex manifold, $\mathcal{A}$ an Azumaya algebra on $X$ and $\alpha$ its class in $Br(X)$. If $E$ is a locally free $\alpha-$twisted coherent sheaf such that $\mathscr{E}nd(E)\simeq\mathcal{A}$, we have an equivalence $$Coh(X,\alpha)\longrightarrow Coh(X,\mathcal{A}),\,\,\,\,\,\,\,\,\,\mathscr{F}\mapsto\mathscr{F}\otimes E^{\vee}.$$
\end{prop}

We now define Chern characters, Mukai vectors and slopes for the objects of $Coh(X,\mathcal{A})$, which allow us to define a notion of stability. For $\mathscr{F}\in Coh(X,\mathcal{A})$ we define $$ch_{\mathcal{A}}(\mathscr{F}):=\frac{ch(\mathscr{F})}{\sqrt{ch(\mathcal{A})}},\,\,\,\,\,\,\,\,v_{\mathcal{A}}(\mathscr{F}):=ch_{\mathcal{A}}(\mathscr{F})\cdot\sqrt{td(X)},$$and if $\omega$ is a K\"ahler class and $\mathscr{F}$ is torsion-free, we let $$\mu_{\mathcal{A},\omega}(\mathscr{F}):=\frac{c_{\mathcal{A},1}(\mathscr{F})\cdot\omega}{rk(\mathscr{F})},$$where $c_{\mathcal{A},1}(\mathscr{F})$ is the component of $ch_{\mathcal{A}}(\mathscr{F})$ in $H^{2}(S,\mathbb{Q})$. We now introduce the notion of stability for $\mathcal{A}-$modules:

\begin{defn}
A torsion-free $\mathscr{F}\in Coh(X,\mathcal{A})$ is $\mu_{\mathcal{A},\omega}-${\em stable} if for every $\mathscr{E}\subseteq\mathscr{F}$ in $Coh(X,\mathcal{A})$ such that $0<rk(\mathscr{E})<rk(\mathscr{F})$, we have $\mu_{\mathcal{A},\omega}(\mathscr{E})<\mu_{\mathcal{A},\omega}(\mathscr{F})$. If $\mu_{\mathcal{A},\omega}(\mathscr{E})\leq\mu_{\mathcal{A},\omega}(\mathscr{F})$ for every such subobject, we say that $\mathscr{F}$ is $\mu_{\mathcal{A},\omega}-${\em semistable}. 
\end{defn}

We notice that if $\mathscr{G}\in Coh(X,\alpha)$ and $E$ is a locally free $\alpha-$twisted sheaf such that $\mathscr{E}nd(E)\simeq\mathcal{A}$, then $ch_{E}(\mathscr{G})=ch_{\mathcal{A}}(\mathscr{G}\otimes E^{\vee})$. It follows that $$v_{E}(\mathscr{G})=v_{\mathcal{A}}(\mathscr{G}\otimes E^{\vee}),\,\,\,\,\,\,\,\mu_{E,\omega}(\mathscr{G})=\mu_{\mathcal{A},\omega}(\mathscr{G}\otimes E^{\vee}),$$so that $\mathscr{F}\in Coh(X,\alpha)$ is $\mu_{E,\omega}-$stable if and only if $\mathscr{F}\otimes E^{\vee}$ is $\mu_{\mathcal{A},\omega}-$stable. 

\begin{oss}
\label{oss:simpson}We notice that $\Lambda:=(\mathscr{O}_{X},\mathcal{A})$ is a sheaf of rings of differential operators following the definition of \cite{Sim}, and $Coh(X,\mathcal{A})$ is the category of $\Lambda-$modules (always in the sense of \cite{Sim}). Moreover, $\mu_{\mathcal{A},\omega}-$stable $\mathcal{A}-$modules are exactly $\mu-$stable $\Lambda-$modules (always in the sense of \cite{Sim}). Even if the definitions of \cite{Sim} are given only for projective manifolds, they can immediately be extended to compact complex manifolds.
\end{oss}

\subsubsection{Twisted sheaves following Yoshioka}

In \cite{Y3} Yoshioka introduces twisted sheaves as a full subcategory of the category of coherent sheaves on a projective bundle. 

More precisely, let $X$ be a complex manifold, $\alpha\in Br(X)$ and $E$ a locally free $\alpha-$twisted sheaf. On an open cover $\mathscr{U}=\{U_{i}\}_{i\in I}$ we represent $E$ by $\{E_{i},\phi_{ij}\}_{i,j\in I}$. Let $\mathbb{P}_{i}:=\mathbb{P}(E_{i})$, together with the map $\pi_{i}:\mathbb{P}_{i}\longrightarrow U_{i}$. The twisted gluing data $\phi_{ij}$ turn to a gluing data $\varphi_{ij}$ of the $\mathbb{P}_{i}$ and of the $\pi_{i}$, getting a projective bundle $\pi:\mathbb{P}\longrightarrow X$ (whose class in $Br(X)$ is $\alpha$). 

As shown in Lemma 1.1 of \cite{Y3}, we have $Ext^{1}(T_{\mathbb{P}/X},\mathscr{O}_{\mathbb{P}})=\mathbb{C}$, hence, up to scalars, there is a unique non-trivial extension $$0\longrightarrow\mathscr{O}_{\mathbb{P}}\longrightarrow G\longrightarrow T_{\mathbb{P}/X}\longrightarrow 0.$$The vector bundle $G$ can be described in another way. Fix a tautological line bundle $\mathscr{O}(\lambda_{i})$ over $\mathbb{P}_{i}$, so that the twisted gluing data $\phi_{ij}$ give isomorphisms $\psi_{ij}:\mathscr{O}(\lambda_{i})\longrightarrow\mathscr{O}(\lambda_{j})$, and $L:=\{\mathscr{O}(\lambda_{i}),\psi_{ij}\}$ is an $\pi^*(\alpha^{-1})-$twisted line bundle on $\mathbb{P}$. Then the vector bundles $\pi_{i}^{*}E_{i}(\lambda_{i})$ glue together to give a locally free sheaf which is isomorphic to $G$.

\begin{defn}
A coherent sheaf $\mathscr{F}$ on $\mathbb{P}$ is called $\mathbb{P}-${\em sheaf} if the canonical morphism $\pi^{*}\pi_{*}(G^{\vee}\otimes\mathscr{F})\longrightarrow G^{\vee}\otimes\mathscr{F}$ is an isomorphism. We let $Coh(\mathbb{P},X)$ be the full subcategory of $Coh(\mathbb{P})$ given by $\mathbb{P}-$sheaves.
\end{defn}

Lemma 1.5 of \cite{Y3} shows that $\mathscr{F}\in Coh(\mathbb{P},X)$ if and only if $\mathscr{F}_{|\mathbb{P}_{i}}\simeq\pi^{*}\mathscr{E}_{|U_i}\otimes\mathscr{O}(\lambda_{i})$ for some $\mathscr{E}\in Coh(U_{i})$. Using this, one shows:

\begin{prop}
\label{prop:equyoshi}Let $X$ be a complex manifold and $\pi:\mathbb{P}\longrightarrow X$ a projective bundle whose class in $Br(X)$ is $\alpha$. Then there is an equivalence of categories $$P:Coh(\mathbb{P},X)\longrightarrow Coh(X,\alpha),\,\,\,\,\,\,\,\,\,P(\mathscr{F}):=\pi_{*}(\mathscr{F}\otimes L^{\vee}).$$
\end{prop}

Following Yoshioka, we have a definition of Chern character, Mukai vector and slope of a $\mathbb{P}-$sheaf $\mathscr{F}$. More precisely, we have 
$$ch_{\mathbb{P}}(\mathscr{F}):=\frac{ch(\pi_{*}(G^{\vee}\otimes\mathscr{F}))}{\sqrt{ch(\pi_{*}(G^{\vee}\otimes G))}},$$so that $$v_{\mathbb{P}}(\mathscr{F})=ch_{\mathbb{P}}(\mathscr{F})\cdot\sqrt{td(S)},\,\,\,\,\,\,\,\mu_{\mathbb{P},\omega}(\mathscr{F}):=\frac{c_{\mathbb{P},1}(\mathscr{F})\cdot\omega}{rk(\mathscr{F})},$$where $c_{\mathbb{P},1}(\mathscr{F})$ is the component of $ch_{\mathbb{P}}(\mathscr{F})$ in $H^{2}(S,\mathbb{Q})$. We now introduce the notion of stability for $\mathbb{P}-$sheaves.

\begin{defn}
We say that a torsion-free $\mathscr{F}\in Coh(\mathbb{P},X)$ is $\mu_{\mathbb{P},\omega}-${\em stable} if for every subobject $\mathscr{E}$ of $\mathscr{F}$ in $Coh(\mathbb{P},X)$ such that $0<rk(\mathscr{E})<rk(\mathscr{F})$, we have $\mu_{\mathbb{P},\omega}(\mathscr{E})<\mu_{\mathbb{P},\omega}(\mathscr{F})$. If $\mu_{\mathbb{P},\omega}(\mathscr{E})\leq\mu_{\mathbb{P},\omega}(\mathscr{F})$ for every such subobject, we say that $\mathscr{F}$ is $\mu_{\mathbb{P},\omega}-${\em semistable}. 
\end{defn}

If $\mathbb{P}=\mathbb{P}(E)$ for some locally free $\alpha-$twisted sheaf $E$, the equivalence $P$ gives $$ch_{\mathbb{P}}(\mathscr{F})=ch_{E}(P(\mathscr{F})),\,\,\,\,\,\,\,v_{\mathbb{P}}(\mathscr{F})=v_{E}(P(\mathscr{F})),$$ $$\mu_{\mathbb{P},\omega}(\mathscr{F})=\mu_{E,\omega}(P(\mathscr{F})).$$It follows that $\mathscr{F}\in Coh(\mathbb{P},X)$ is $\mu_{\mathbb{P},\omega}-$stable if and only if $P(\mathscr{F})$ is $\mu_{E,\omega}-$stable. 

If $\mathscr{F}$ is a $\mu_{\mathbb{P},\omega}-$stable $\mathbb{P}-$sheaf, as $Coh(\mathbb{P},S)$ is a full subcategory of $Coh(\mathbb{P})$ and as the functor $P$ is an equivalence, we have that $$Ext^{1}_{Coh(\mathbb{P},S)}(\mathscr{F},\mathscr{F})\simeq Ext^{1}_{Coh(S,\alpha)}(P(\mathscr{F}),P(\mathscr{F})),$$and $$Ext^{2}_{Coh(\mathbb{P},S)}(\mathscr{F},\mathscr{F})\simeq Ext^{2}_{Coh(S,\alpha)}(P(\mathscr{F}),P(\mathscr{F})).$$

\subsubsection{Chern classes following Huybrechts and Stellari}

If we consider twisted sheaves following C\u ald\u araru, there is another possible definition of their Chern classes and character, introduced by Huybrechts and Stellari in \cite{HS2}, that we recall here.

Consider a complex manifold $X$ and $\alpha\in H^{2}(X,\mathcal{O}^{*}_{X})$, and fix a \v Cech $2-$cocycle $\{\alpha_{ijk}\}$ representing $\alpha$ on an open covering $\{U_{i}\}_{i\in I}$ of $X$. Moreover, choose a \v Cech $2-$cocyle $\{B_{ijk}\}$, where $B_{ijk}\in\Gamma(U_{ijk},\mathbb{Q})$, such that $\alpha_{ijk}=\exp(B_{ijk})$ (viewed as local sections of $\mathbb{R}/\mathbb{Z}=U(1)\subseteq\mathcal{O}_{X}^{*}$).

As the sheaf $\mathcal{C}^{\infty}_{X}$ of $C^{\infty}-$functions on $X$ is acyclic, up to supposing the covering $\{U_{i}\}_{i\in I}$ is sufficiently fine, there are $a_{ij}\in\Gamma(U_{ij},C^{\infty})$ such that $$B_{ijk}=-a_{ij}+a_{ik}-a_{jk}.$$

Now, let us consider an $\alpha-$twisted sheaf given by $\mathscr{F}=\{\mathscr{F}_{i},\phi_{ij}\}$ and let $$\psi_{ij}:=\phi_{ij}\cdot\exp(a_{ij}),$$which is clearly an isomorphism between the restrictions of $\mathscr{F}_{i}$ and $\mathscr{F}_{j}$ to $U_{ij}$. It is moreover easy to show that $$\psi_{ij}\circ\psi_{jk}\circ\psi_{ki}=id,$$hence the sheaf $\mathscr{F}_{B}=\{\mathscr{F}_{i},\psi_{ij}\}_{i,j\in I}$ is an untwisted sheaf. We then let $$ch^{B}(\mathscr{F}):=\ch(\mathscr{F}_{B}).$$The definition given in this way depends only on $B$. 

The relation between $ch^{B}$ and the previous Chern characters is explained in \cite{HS}, and goes as follows, supposing that $\alpha\in Br(X)$. Let $E$ be a locally free $\alpha-$twisted sheaf and $$B_{E}:=\frac{c_{1}^{B}(E)}{rk(E)},$$where $c_{1}^{B}(E)$ is the degree two part of $ch^{B}(E)$. Then we have $$ch^{B}(\mathscr{F})=ch_{E}(\mathscr{F})\cdot\exp(B_{E}).$$

\subsection{Genericity for polarizations}

We now extend the notion of genericity for polarization to the twisted case. As we did in section 2.2, we first introduce a notion of discriminant for twisted sheaves, which depends on the choice of a locally free $E\in Coh(S,\alpha)$. 

\subsubsection{Discriminant of a twisted sheaf}

If $\mathscr{F}$ is an $\alpha-$twisted coherent sheaf, we call \textit{discriminant} of $\mathscr{F}$ with respect to $E$ the number $$\Delta_{E}(\mathscr{F}):=\frac{1}{2rk^{2}(\mathscr{F})}v_{E}^{2}(\mathscr{F})+1.$$If $\alpha=0$, this is just $\Delta(\mathscr{F})$ by equation (\ref{eq:v2}). More generally, the discriminant does not depend on $E$, as shown in the following:

\begin{lem}
\label{lem:nondip}
Let $\alpha\in Br(S)$ and $\mathscr{F}\in Coh(S,\alpha)$. If $E_{1},E_{2}\in Coh(S,\alpha)$ are locally free, then $\Delta_{E_{1}}(\mathscr{F})=\Delta_{E_{2}}(\mathscr{F})$.
\end{lem}

\proof Let $E\in Coh(S,\alpha)$ be locally free of rank $s$, and pose $r:=rk(\mathscr{F})$, $\xi:=c_{1}(\mathscr{F}\otimes E^{\vee})$, $a:=ch_{2}(\mathscr{F}\otimes E^{\vee})$ and $b:=ch_{2}(E\otimes E^{\vee})$. By equation (\ref{eq:vsquare}) we have $$\Delta_{E}(\mathscr{F})=\frac{1}{2r^{2}}\bigg(\frac{\xi^{2}}{s^{2}}-\frac{2ra}{s}+\frac{r^{2}b}{s^{2}}-2r^{2}\bigg)+1.$$An easy computation shows that $$\frac{\xi^{2}}{s}-\frac{2ra}{s}+\frac{r^{2}b}{s^{2}}=-\frac{ch_{2}(\mathscr{F}\otimes\mathscr{F}^{\vee}\otimes E\otimes E^{\vee})}{s^{2}}+\frac{r^{2}ch_{2}(E\otimes E^{\vee})}{2s^{2}}=$$ $$=-ch_{2}(\mathscr{F}\otimes\mathscr{F}^{\vee}),$$so that 
\begin{equation}
\label{eq:deltae}
\Delta_{E}(\mathscr{F})=\frac{1}{2r^{2}}(-ch_{2}(\mathscr{F}\otimes\mathscr{F}^{\vee})-2r^{2})+1,
\end{equation}
which does not depend on $E$, implying the statement.\endproof

For $v\in H^{2*}(S,\mathbb{Q})$ and $\alpha\in Br(S)$, we let $$\Delta_{\alpha}(v):=\Delta_{\mathscr{F}^{\vee\vee}}(\mathscr{F}),$$where $\mathscr{F}\in Coh(S,\alpha)$ is torsion free and $v_{\alpha}(\mathscr{F})=v$. By Lemma \ref{lem:nondip} this is well defined and if $\alpha=0$, then $\Delta_{0}(v)=\Delta(v)$. We now prove a generalization to twisted sheaves of the Bogomolov inequality:

\begin{prop}
\label{prop:deltatwisted}Let $\alpha\in Br(S)$, $\mathscr{F}\in Coh(S,\alpha)$ and $\omega$ a K\"ahler class on $S$. If $\mathscr{F}$ is $\mu_{\alpha,\omega}-$semistable, then $\Delta_{\alpha}(\mathscr{F})\geq 0$.
\end{prop}

\proof It is easy to see that $\mathscr{F}$ is $\mu_{\alpha,\omega}-$semistable if and only if $\mathscr{F}^{\vee}$ is $\mu_{\alpha^{-1},\omega}-$semistable. In particular, this implies that $\mathscr{F}$ is $\mu_{\alpha,\omega}-$semistable if and only if $\mathscr{F}^{\vee\vee}$ is $\mu_{\alpha,\omega}-$semistable. 

Now, notice that $\mathscr{F}^{\vee\vee}\otimes\mathscr{F}^{\vee}=(\mathscr{F}\otimes\mathscr{F}^{\vee})^{\vee\vee}$, hence if $l$ is the length of the singular locus of $\mathscr{F}\otimes\mathscr{F}^{\vee}$, it follows that $$ch_{2}(\mathscr{F}^{\vee\vee}\otimes \mathscr{F}^{\vee})=ch_{2}(\mathscr{F}\otimes\mathscr{F}^{\vee})+l.$$By equation (\ref{eq:vsquare}), it then follows that $$v_{\alpha}^{2}(\mathscr{F}^{\vee\vee})=v_{\alpha}^{2}(\mathscr{F})-2l\leq v_{\alpha}^{2}(\mathscr{F}).$$As $rk(\mathscr{F})=rk(\mathscr{F}^{\vee\vee})$ it follows that $\Delta_{\alpha}(\mathscr{F})\geq \Delta_{\alpha}(\mathscr{F}^{\vee\vee})$, hence we just need to show the statement for $\mathscr{F}^{\vee\vee}$.

Let now $F:=\mathscr{F}^{\vee\vee}$, which is locally free and $\mu_{\alpha,\omega}-$semistable. By the Kobayashi-Hitchin correspondence for twisted sheaves as proved by Wang in \cite{Wang}, the sheaf $\mathscr{E}nd(F)=F\otimes F^{\vee}$ is $\mu_{\omega}-$polystable, so that $\Delta(F\otimes F^{\vee})\geq 0$ by the Bogomolov inequality.

Choose now a locally free $E\in Coh(S,\alpha)$ of rank $s$, and let $b:=ch_{2}(E\otimes E^{\vee})$. By Lemma \ref{lem:nondip} we have $\Delta(F\otimes F^{\vee})=\Delta_{E\otimes E^{\vee}}(F\otimes F^{\vee})$, so that $\Delta_{E\otimes E^{\vee}}(F\otimes F^{\vee})\geq 0$. If $\xi=c_{1}(F\otimes E^{\vee})$ and $a=ch_{2}(F\otimes E^{\vee})$, it follows from equation (\ref{eq:vsquare}) that $$v_{E\otimes E^{\vee}}^{2}(F\otimes F^{\vee})=\frac{2r^{2}\xi^{2}}{s^{2}}-\frac{4r^{3}a}{s}+\frac{2r^{2}b}{s^{2}}-2r^{4}=2r^{2}v_{E}^{2}(F)+2r^{4}.$$Hence $$0\leq\Delta_{E\otimes E^{\vee}}(F\otimes F^{\vee})=\frac{1}{2r^{4}}v^{2}_{E\otimes E^{\vee}}(F\otimes F^{\vee})+1=2\Delta_{E}(F).$$Hence $\Delta_{\alpha}(F)=\Delta_{E}(F)\geq 0$, and we are done.\endproof

\subsubsection{Walls and chambers}

Now, let $$W_{\alpha,v}:=\{D\in NS(S)\,|\,-\frac{r^{4}}{2}\Delta_{\alpha}(v)\leq D^{2}<0\}.$$If $\alpha=0$, we have $W_{0,v}=W_{v}$. 

\begin{defn}
If $D\in W_{\alpha,v}$, we call the hyperplane $D^{\perp}$ an $(\alpha,v)-${\em wall}. A connected component of $\mathscr{K}_{S}\setminus\bigcup_{D\in W_{\alpha,v}}D^{\perp}$ is called $(\alpha,v)-${\em chamber}. A polarization $\omega\in\mathscr{K}_{S}$ is $(\alpha,v)-${\em generic} if it lies in a $(\alpha,v)-$chamber.
\end{defn}

A polarization $\omega$ is $(0,v)-$generic if and only if it is $v-$generic. We are now ready to prove one of the main results of this section about changing polarization inside a chamber. The argument is the same one for untwisted sheaves, here adapted to the twisted case.

\begin{prop}
\label{prop:gentwisted}Let $\alpha\in Br(S)$, $v\in H^{2*}(S,\mathbb{Q})$ and $\omega,\omega'$ two $(\alpha,v)-$generic polarizations lying in the same $(\alpha,v)-$chamber. If $\mathscr{F}\in Coh(S,\alpha)$ is a torsion free sheaf such that $v_{\alpha}(\mathscr{F})=v$, then $\mathscr{F}$ is $\mu_{\alpha,\omega}-$stable if and only if it is $\mu_{\alpha,\omega'}-$stable.
\end{prop}

\proof The proof is divided in two steps.

\textit{Step 1}. Choose an $\alpha-$twisted locally free sheaf $E$, and let $r:=rk(\mathscr{F})$, $\xi:=c_{1}(\mathscr{F}\otimes E^{\vee})$, $a:=ch_{2}(\mathscr{F}\otimes E^{\vee})$, $s:=rk(E)$ and $b:=ch_{2}(E\otimes E^{\vee})$. Let $\omega$ be any polarization, and suppose that $\mathscr{F}$ is properly $\mu_{\alpha,\omega}-$semistable: hence there is an $\alpha-$twisted subsheaf $\mathscr{E}\subseteq\mathscr{F}$ such that $\mu_{E,\omega}(\mathscr{E})=\mu_{E,\omega}(\mathscr{F})$. We let $r':=rk(\mathscr{E})$, $\xi':=c_{1}(\mathscr{E}\otimes E^{\vee})$ and $a':=ch_{2}(\mathscr{E}\otimes E^{\vee})$, where $0<r'<r$ and $\xi'\in NS(S)\otimes\mathbb{Q}$. Moreover, let $$D:=r\frac{\xi'}{s}-r'\frac{\xi}{s},$$so that $D\cdot\omega=0$. Hence $D^{2}\leq 0$, as $\omega$ is a K\"ahler form. 

Now, let $\mathscr{G}:=\mathscr{F}/\mathscr{E}$, and we suppose without loss of generality that $\mathscr{E}$ is saturated and that $\mathscr{E}$ and $\mathscr{G}$ are $\mu_{\alpha,\omega}$-semistable. Moreover, let $r'':=rk(\mathscr{G})$, $\xi'':=c_{1}(\mathscr{G}\otimes E^{\vee})$ and $a'':=ch_{2}(\mathscr{G}\otimes E^{\vee})$, so that $r''=r-r'$, $\xi''=\xi-\xi'$ and $a''=a-a'$. Finally, let $$K:=\frac{v_{\alpha}^{2}(\mathscr{F})}{r}-\frac{v_{\alpha}^{2}(\mathscr{E})}{r'}-\frac{v_{\alpha}^{2}(\mathscr{G})}{r''}.$$

We notice that as $\mathscr{E}$ and $\mathscr{G}$ are $\mu_{\alpha,\omega}-$semistable, by Proposition \ref{prop:deltatwisted} we have $\Delta_{\alpha}(\mathscr{E}),\Delta_{\alpha}(\mathscr{G})\geq 0$, meaning $v_{\alpha}^{2}(\mathscr{E})\geq -2(r')^{2}$ and $v_{\alpha}^{2}(\mathscr{G})\geq -2(r'')^{2}$. Hence we get 
\begin{equation}
\label{eq:ineqk}
K\leq\frac{v_{\alpha}^{2}(\mathscr{F})}{r}+2r.
\end{equation} 
On the other hand, by equation (\ref{eq:vsquare}) we have $$K=\frac{\xi^{2}}{rs^{2}}-\frac{(\xi')^{2}}{r's^{2}}-\frac{(\xi'')^{2}}{r''s^{2}}=-\frac{r^{2}(\xi')^{2}+(r')^{2}\xi^{2}-2rr'\xi\xi'}{s^{2}rr'r''}.$$By definition of $D$ we have $$D^{2}=\frac{r^{2}(\xi')^{2}+(r')^{2}\xi^{2}-2rr'\xi\xi'}{s^{2}},$$so that the inequality (\ref{eq:ineqk}) implies $$D^{2}=-rr'r''K\geq-r'(r-r')v_{\alpha}^{2}(\mathscr{F})-2r^{2}r'(r-r').$$But as $r'(r-r')\leq r^{2}/4$, we finally get $$D^{2}\geq-\frac{r^{2}}{4}v_{\alpha}^{2}(\mathscr{F})-\frac{r^{4}}{2}=-\frac{r^{4}}{2}\Delta_{\alpha}(\mathscr{F})=-\frac{r^{4}}{2}\Delta_{\alpha}(v).$$In conclusion, $D\in W_{\alpha,v}\cup\{0\}$.

\textit{Step 2}. Suppose that $\mathscr{F}$ is $\mu_{\alpha,\omega}-$stable but not $\mu_{\alpha,\omega'}-$stable. Let $$[\omega,\omega']:=\{\omega_{t}:=t\omega'+(1-t)\omega\,|\,t\in[0,1]\}$$be the segment from $\omega$ to $\omega'$, and let $B_{\alpha}$ be the family of subsheaves of $\mathscr{F}$ whose slope with respect to $E$ and $\omega'$ is bounded from below. 

If $\mathscr{E}\in B_{\alpha}$, then $\mathscr{E}\otimes E^{\vee}$ is a subsheaf of $\mathscr{F}\otimes E^{\vee}$, and $\mu_{E,\omega}(\mathscr{E})=\mu_{\omega}(\mathscr{E}\otimes E^{\vee})$. This implies that $\mathscr{E}\otimes E^{\vee}$ is in the family $B$ of subsheaves of $\mathscr{F}\otimes E^{\vee}$ whose slope with respect to $\omega$ is bounded from below. As the family $B$ is bounded, it follows that the family $B_{\alpha}$ is bounded. Using the same argument as in the proof of Lemma \ref{lem:openness}, one can then conclude that there is $t\in]0,1]$ such that $\mathscr{F}$ is properly $\mu_{\alpha,\omega_{t}}-$semistable.

Hence there is a subsheaf $\mathscr{E}$ of $\mathscr{F}$ of rank $0<s<r$ such that $\mu_{E,\omega_{t}}(\mathscr{E})=\mu_{E,\omega_{t}}(\mathscr{F})$. If $D=rc_{E,1}(\mathscr{E})-sc_{E,1}(\mathscr{F})$, it follows that $D\cdot\omega_{t}=0$. As $D\cdot\omega\neq 0$, we have $D\neq 0$, hence $D^{2}<0$. But as $\mathscr{F}$ is $\mu_{E,\omega_{t}}-$semistable, Step 1 implies that $D\in W_{\alpha,v}$, which is not possible as $\omega_{t}$ is in the same $(\alpha,v)-$chamber as $\omega$ and $\omega'$. In conclusion, the sheaf $\mathscr{F}$ has to be $\mu_{E,\omega'}-$stable, and we are done.\endproof 

\subsection{Moduli space of stable twisted sheaves}

We now introduce (relative) moduli spaces of stable twisted sheaves. On projective manifolds these were constructed by Yoshioka in \cite{Y3} (viewing twisted sheaves as $\mathbb{P}-$sheaves, and using a GIT construction), and by Lieblich in \cite{L} (viewing twisted sheaves as sheaves on some $\mathscr{O}^{*}-$gerbe).

Here we first provide a relative moduli space of simple twisted sheaves by viewing them as simple $\mathbb{P}-$sheaves. The relative moduli space of stable sheaves will then be an open subset of it.

\subsubsection{The relative moduli space of simple twisted sheaves}

Consider a smooth and proper morphism $\pi:\mathscr{X}\longrightarrow T$ such that for every $t\in T$ the fibre $X_{t}$ over $t$ is a K3 surface. We assume for simplicity that $T$ is a complex manifold, although the constructions work over complex spaces as well.

Suppose moreover that we are given a complex manifold $\mathscr{P}$ together with a morphism $f:\mathscr{P}\longrightarrow\mathscr{X}$ of $T-$complex spaces such that for every $t\in T$, the morphism $f_{t}:P_{t}\longrightarrow X_{t}$ is a projective bundle, where $P_{t}=f^{-1}(X_{t})$. For every $t\in T$ the projective bundle $P_{t}\longrightarrow X_{t}$ defines a class $\alpha_{t}$ in the Brauer group $Br(X_{t})$.

Now, let $f':=\pi\circ f$, so that we get a map $f':\mathscr{P}\longrightarrow T$. By Theorem (6.4) of \cite{KO}, there is a complex space $\mathscr{M}(\mathscr{P}/T)$ together with a holomorphic surjective map $$\phi:\mathscr{M}(\mathscr{P}/T)\longrightarrow T$$which is a relative moduli space of simple coherent sheaves on $\mathscr{P}$: for every $t\in T$ the fibre $\mathscr{M}_{t}$ of $\phi$ over $t$ is the moduli space of simple coherent sheaves on $P_{t}$.

Now, $\mathscr{F}\in Coh(P_{t})$ is simple if and only if $End(\mathscr{F})\simeq\mathbb{C}$. As $Coh(P_{t},X_{t})$ is a full subcategory of $Coh(P_{t})$, a $P_{t}-$sheaf $\mathscr{F}$ is simple in $Coh(P_{t},X_{t})$ if and only if it is simple in $Coh(P_{t})$. Hence simple $P_{t}-$sheaves form a subset $\mathscr{M}^{s}(\mathscr{P}/T)$ of $\mathscr{M}(\mathscr{P}/T)$. 

As the condition defining $\mathbb{P}-$sheaves is open (see Lemma 1.5 of \cite{Y3}), it follows that $\mathscr{M}^{s}(\mathscr{P}/T)$ is open in $\mathscr{M}(\mathscr{P}/T)$, hence it is a complex space together with a holomorphic map $\psi:\mathscr{M}^{s}(\mathscr{P}/T)\longrightarrow T$. This is the relative moduli space of simple $\mathscr{P}-$sheaves on $\mathscr{X}$. 

The relative projective bundle $f:\mathscr{P}\longrightarrow\mathscr{X}$ corresponds to the existence of a relative Azumaya algebra $\mathcal{A}$ on $\mathscr{X}$: for every $t\in T$, we have $P_{t}=\mathbb{P}(E_{t})$ for some locally free $\alpha_{t}-$twisted sheaf on $X_{t}$, and we let $\mathcal{A}_{t}:=E_{t}\otimes E_{t}^{\vee}$. The previous equivalence of categories of twisted sheaves then shows that $\mathscr{M}^{s}(\mathscr{P}/T)$ is the relative moduli space of simple $\mathcal{A}-$modules on $\mathscr{X}$ or, equivalently, the relative moduli space of simple twisted sheaves on $\mathscr{X}$.

\subsubsection{The relative moduli space of stable twisted sheaves}\label{mod-rel-stab}

We now produce out of $\psi:\mathscr{M}^{s}(\mathscr{P}/T)\longrightarrow T$ the relative moduli space of stable twisted sheaves. Choose $v=(v_{0},v_{1},v_{2})\in H^{2*}(S,\mathbb{Q})$ such that $v_{1}\in NS(S_{t})$ for every $t\in T$, and $v_{0}\geq 2$. We let $\mathscr{M}^{s}_{v}(\mathscr{P}/T)$ be the component of $\mathscr{M}^{s}(\mathscr{P}/T)$ parametrizing simple $\mathbb{P}-$sheaves of Mukai vector $v$, and we write $\psi_{v}:\mathscr{M}^{s}_{v}(\mathscr{P}/T)\longrightarrow T$ for $\psi_{|\mathscr{M}^{s}_{v}(\mathscr{P}/T)}$.

In order to define the moduli space of stable twisted sheaves of Mukai vector $v$, we need a section $\widetilde{\omega}\in R^{2}\pi_{*}\mathbb{C}$ such that $\omega_{t}:=\widetilde{\omega}_{|X_{t}}$ is a K\"ahler class for every $t\in T$, which is used to define stability on every fibre. As stable twisted sheaves are simple, we let $\mathscr{M}^{\mu}_{v}(\mathscr{P}/T,\widetilde{\omega})$ be the subset of $\mathscr{M}^{s}_{v}(\mathscr{P}/T)$ whose fibre over $t\in T$ is given by the simple $P_{t}-$sheaves which are $\mu_{P_{t},\omega_{t}}-$stable and whose Mukai vector is $v$. We then have a natural map (of sets) $$p:\mathscr{M}^{\mu}_{v}(\mathscr{P}/T,\widetilde{\omega})\longrightarrow T.$$ 

The main result of this section is the following

\begin{prop}
\label{prop:modstable}Let $\pi:\mathscr{X}\longrightarrow T$, $f:\mathscr{P}\longrightarrow\mathscr{X}$, $v$ and $\widetilde{\omega}$ be as before. Then $\mathscr{M}^{\mu}_{v}(\mathscr{P}/T,\widetilde{\omega})$ is an open subset of $\mathscr{M}^{s}_{v}(\mathscr{P}/T)$. Hence it is a complex manifold, and the map $p$ is holomorphic.
\end{prop}

\proof By Remark \ref{oss:simpson}, the openness can be proved as in Lemma 3.7 of \cite{Sim}. Indeed, if $\mathscr{F}\in Coh(\mathbb{P},S)$ and $F:=P(\mathscr{F})^{\vee\vee}$, then $\mathscr{F}$ is $\mu_{\mathbb{P},\omega}-$stable if and only if $P(\mathscr{F})\otimes P(F)^{\vee}$ is $\mu_{\mathcal{A},\omega}-$stable in $Coh(S,\mathcal{A})$, where $\mathcal{A}=P(F)\otimes P(F)^{\vee}$. Moreover, the openness of stability in the analytic case may be proved in the usual way, by using boundedness results which are contained in \cite{To2} and \cite{To-inprep}. The separatedness follows from Proposition (6.6) of \cite{KO} since the parameterized sheaves are stable.\endproof

Standard deformation arguments following \cite{BF} allow us to show that if $p:\mathscr{M}:=\mathscr{M}^{\mu}_{v}(\mathscr{P}/\mathscr{X},\widetilde{\omega})\longrightarrow T$ is the relative moduli space of twisted stable sheaves, then for every $t\in T$ and for every $\mathscr{F}\in p^{-1}(t)=\mathscr{M}_{t}$ we have $$T_{[\mathscr{F}]}\mathscr{M}_{t}\simeq Ext^{1}_{Coh(X_{t},\alpha_{t})}(P(\mathscr{F}),P(\mathscr{F})),$$that the obstruction for the existence of deformations of $\mathscr{F}$ live in $$Ext^{2}_{Coh(X_{t},\alpha_{t})}(P(\mathscr{F}),P(\mathscr{F})),$$and that we have an exact sequence 
\begin{equation}
\label{eq:extan}
Ext^{1}_{Coh(X_{t},\alpha_{t})}(P(\mathscr{F}),P(\mathscr{F}))\longrightarrow T_{[\mathscr{F}]}\mathscr{M}\longrightarrow
\end{equation}
\begin{equation*}
\longrightarrow T_{t}T\longrightarrow Ext^{2}_{Coh(X_{t},\alpha_{t})}(P(\mathscr{F}),P(\mathscr{F}))_{0}.
\end{equation*}
It follows from this exact sequence and by the previous discussion, that the morphism $p:\mathscr{M}\longrightarrow T$ is smooth.
 
If $T$ is reduced to a point, then $X$ is just a K3 surface $S$ and $P\longrightarrow S$ is a projective bundle whose class in $Br(S)$ is $\alpha$. The \textit{moduli space of $\mu_{\alpha,\omega}-$stable $\alpha-$twisted sheaves of twisted Mukai vector $v$ on $S$} will then be denoted $M^{\mu}_{\alpha,v}(S,\omega)$.

\begin{oss}
\label{oss:isotwist}
Suppose that $\alpha=0$ and let $$\gamma:=\frac{ch(\mathscr{F}^{\vee})}{\sqrt{ch(\mathscr{F}\otimes\mathscr{F}^{\vee})}}$$for $\mathscr{F}\in M^{\mu}_{0,v}(S,\omega)$. Then $v_{0}(\mathscr{F})=v$ if and only if $v(\mathscr{F})=v/\gamma$, so that $M^{\mu}_{0,v}(S,\omega)\simeq M^{\mu}_{v/\gamma}(S,\omega)$. We even notice that $\omega$ is $(0,v)-$generic if and only if it is $v/\gamma-$generic.
\end{oss}

\subsubsection{Moduli spaces of stable twisted sheaves over projective K3 surfaces}

If the base K3 surface $S$ is projective, from \cite{Y3} we have some informations more about the moduli spaces of stable twisted sheaves. We make use of the following notation: let $\alpha\in Br(S)$ and $\mathscr{F}$ a torsion free $\alpha-$twisted sheaf whose twisted Mukai vector is $w=(r,0,a)$. 

We let $F$ be a locally free $\alpha-$twisted sheaf and $\xi$ be a representative of the class of $\mathbb{P}(E)$ in $H^{2}(S,\mathbb{Z})$. We let $e^{\xi/r}:=(1,\xi/r,\xi^{2}/2r^{2})$ and $w_{\xi}:=e^{\xi/r}\cdot w$, so that $$w_{\xi}=(r,\xi,a+\xi^{2}/2r).$$It is worthwhile to notice that there is a topological vector bundle $E_{\xi}$ on $S$ such that $v(E_{\xi})=w_{\xi}$. As shown in \cite{Y3}, we have $w_{\xi}\in H^{2}(S,\mathbb{Z})$ (while in general we have $w\in H^{2}(S,\mathbb{Q})$). 

\begin{oss}
\label{oss:wxi}
If $\alpha=0$ and $\mathscr{F}$ is a $\mu_{\omega}-$stable sheaf whose Mukai vector is $v=(r,\xi,a)$, write $a=c+r$ where $c=ch_{2}(\mathscr{F})$. The $0-$twisted Mukai vector of $\mathscr{F}$ is then $w=(r,0,r+a'/2r)$, where $a'=ch_{2}(\mathscr{F}\otimes\mathscr{F}^{\vee\vee})$. We notice that $a'=2rc-\xi^{2}$, hence $w=(r,0,r+c-\xi^{2}/2r)$. A representative of the class of $\mathbb{P}(E)$ in this case can be chosen to be $\xi$ itself. Hence we have $$w_{\xi}=e^{\xi/r}w=(r,\xi,r+c)=v.$$
\end{oss}

The following is Theorem 3.16 of \cite{Y3}:

\begin{thm}
\label{thm:modtwiamp}Let $S$ be a projective K3 surface, $w=(r,\zeta,b)\in H^{2}(S,\mathbb{Q})$ and $\alpha\in Br(S)$. Choose a representative $\xi$ of $\alpha$ in $H^{2}(S,\mathbb{Z})$, and suppose that $w_{\xi}$ is primitive. Moreover, let $H$ be a $(\alpha,w)-$generic ample line bundle on $S$. Then the moduli space $M^{s}_{\alpha,w}(S,H)$ is an irreducible symplectic manifold which is deformation equivalent to a Hilbert scheme of points on $S$.
\end{thm}

We have the following result, which is the twisted version of Theorem \ref{thm:change}:

\begin{prop}
\label{prop:modtwiproj}Let $S$ be a projective K3 surface, $w=(r,\zeta,b)$ a Mukai vector and $\alpha\in Br(S)$. Choose $\xi$ to be a representative of $\alpha$ in $H^{2}(S,\mathbb{Z})$, and suppose that $r$ and $\xi$ are prime to each other. If $\omega$ is a $(\alpha,w)-$generic polarization, then $M^{\mu}_{\alpha,w}(S,\omega)$ is an irreducible symplectic manifold which is deformation equivalent to a Hilbert scheme of points on $S$.
\end{prop}

\proof By Lemma \ref{lem:projomega}, Proposition \ref{prop:gentwisted} and following the same proof of Theorem \ref{thm:change}, we see that there is an ample line bundle $H$ such that $M^{\mu}_{\alpha,w}(S,\omega)=M^{\mu}_{\alpha,w}(S,H)$. This last moduli space is an irreducible symplectic manifold which is deformation equivalent to a Hilbert scheme of points on $S$ by Theorem \ref{thm:modtwiamp}.\endproof

\subsubsection{Quasi-universal families}

We conclude this section with the following result about the existence of a quasi-universal family; cf. \cite{AMT} for the absolute untwisted case.

\begin{prop}
\label{prop:univfam}
Let $\pi:\mathscr{X}\longrightarrow T$, $f:\mathscr{P}\longrightarrow\mathscr{X}$, $v=(v_{0},v_{1},v_{2})$ and $\widetilde{\omega}$ be as before. Let $\mathcal{A}$ be a relative Azumaya algebra corresponding to $\mathscr{P}$, and for every $t\in T$ let $\alpha_{t}\in Br(X_{t})$ be the class of $\mathcal{A}_{t}$. Suppose that there is a locally free $\mathcal{A}$-module $\mathcal{V}$ verifying the two following properties for every $t\in T$:
\begin{enumerate}
 \item the restriction $\mathcal{V}_{t}$ of $\mathcal{V}$ to $X_{t}$ is $\mu_{\alpha_{t},\omega_{t}}-$stable;
 \item the twisted Mukai vector of $\mathcal{V}_{t}$ is $(v_{0},v_{1},w_{2})$, where $w_{2}<v_{2}$.
\end{enumerate}
Then there is a quasi-universal family on $\mathscr{M}^{\mu}_{v}(\mathscr{P}/T,\widetilde{\omega})\times_{T}\mathscr{X}$.
\end{prop}

\proof Let $\mathscr{M}:=\mathscr{M}^{\mu}_{\widetilde{v}}(\mathscr{P}/T,\widetilde{\omega})$. 
As for stable coherent sheaves, there is an open covering $\mathscr{U}=\{U_{i}\}_{i\in I}$ of $\mathscr{M}$ given by analytic subsets endowed with universal $\mathcal{A}$-modules $\mathscr{F}_{i}$.

Let $p_i:U_i\times_T\mathscr{X}\to U_i$ and $q_i:U_i\times_T\mathscr{X}\to \mathscr{X}$ denote the two projections.
We  put $\mathscr{E}_{i}:=\mathscr{F}_{i}\otimes_{q_i^*\mathcal{A}} q_{i}^{*}\mathcal{V}^{\vee}$. 
By the choice of $\mathcal{V}$ we have $R^0p_{i,*}\mathscr{E}_i=0=R^2p_{i,*}\mathscr{E}_i$ and $W_{i}:=R^1p_{i,*}\mathscr{E}_i$ is a non-trivial locally free $\mathscr{O}_{U_i}$-module whose rank is independent of $i$.

It is easy to check now that the $\mathcal{A}$-modules $\mathscr{F}_i\otimes_{\mathscr{O}}p_i^*W_i^{\vee}$ glue together to give the desired quasi-universal family.
\endproof

\subsection{Deformation of stable twisted sheaves along twistor lines}

In this subsection we describe and generalize a construction used by several authors in the case of stable locally free sheaves of slope zero, cf. \cite{To99}, \cite{Ver}, \cite{Ver08}, \cite{Mar}.

Let $(S, I, \omega)$ be a polarized K3 surface and $\pi:Z(S)\longrightarrow\mathbb{P}^{1}$ its twistor family. We suppose that the fibre over 0 is $S_{0}=(S,I)$, and we write $S_{t}=(S,I_{t})$ for the fibre over $t$. Here $I=I_0$ and $I_t$ denote the complex structures on $S$. With this convention we have $S_\infty=(S, I_\infty)=(S, -I)$. Recall that the choice of $\omega$ on $(S,I)$ is equivalent to the choice of a Riemannian metric $g$ which is compatible with $I$ and whose associated K\"ahler class is $\omega$. Along the twistor line the metric $g$ remains compatible with $I_{t}$, the associated class $\omega_{t}$ is K\"ahler, and we get a section $\widetilde{\omega}$ of $R^{2}\pi_{*}\mathbb{C}$ which is $\omega_{t}$ on $S_{t}$. Slope stability on $S_t$ will be considered with respect to $\omega_t$. 

Before we introduce deformations of sheaves along twistor lines we make an observation on $(1,1)$-forms on the twistor space of $S$. Recall that, as a differentiable manifold, $Z(S)$ is the product $S\times \P^1$, which is endowed with a complex structure in the following way (see \cite{HKLR}): cover $\P^1$ by two charts (each isomorphic to $\C$) and take $\zeta$ the complex coordinate function on one of them and $\zeta^{-1}$ on the other. Further, let $I,J,K$ be the complex structures on $S$ which make it into a hyperk\"ahler manifold. If $I_{\P^1}$ is the complex structure on $\P^1$ then put the following complex structure to act on the tangent space $T_{S}\times T_{\mathbb{P}^{1}}$ of $S\times \P^1$:

$$\mathfrak{I}:=\bigg(\frac{1-\zeta\bar\zeta}{1+\zeta\bar\zeta}I+ \frac{\zeta+\bar\zeta}{1+\zeta\bar\zeta}J+   i\frac{\zeta-\bar\zeta}{1+\zeta\bar\zeta}K, I_{\P^1}\bigg). $$

With respect to this complex structure the projection $q: S\times \P^1\to S$ is not holomorphic but only $C^\infty$.

\begin{lem}
\label{asd} 
Let $\psi$ be a $(1,1)$-form on $(S, I,\omega)$. Its pull-back $q^*\psi$ is a $(1,1)$-form on $Z(S)$ if and only if $\psi$ is anti-self-dual on $(S, I,\omega)$. 
\end{lem}

\proof Let $\Psi:=q^*\psi$. It is a $2$-form on $Z(S)$, so it is of type $(1,1)$ if and only if
$\Psi(\mathfrak{I}v,\mathfrak{I}w)=\Psi(v,w)$ for any pair $(v,w)$ of real tangent vectors at a point of $Z(S)$. As $\mathfrak{I}$ preserves the horizontal and the vertical directions on $Z(S)=S\times \P^1$, and as $\Psi(v,w)=0$ if one of the tangent vectors $v$ or $w$ is horizontal, it suffices to check $\Psi(\mathfrak{I}v,\mathfrak{I}w)=\Psi(v,w)$ only on vertical vectors, meaning that the restrictions of $\Psi$ to the fibres of $\pi:Z(S)\longrightarrow\mathbb{P}^{1}$ are of type $(1,1)$.

Suppose that $\psi$ is anti-self-dual. This property only depends on $g$ and on the orientation of $S$: as $g$ is compatible with each complex structure $I_t$, it follows that the restriction of $\Psi$ to each fibre of $\pi$ is then anti-self-dual. In particular, it is of type $(1,1)$, hence also $\Psi$ is of type $(1,1)$ on $Z(S)$. 

Conversely, if $\psi$ is not anti-self-dual, then it decomposes as $\psi=\psi^{SD}+\psi^{ASD}$, where the self-dual part is of the form $\psi^{SD}=f\omega_I$ for some non-zero function $f$. But $\omega_I $ is not of type $(1,1)$ with respect to $J$ so neither will be $\Psi$.
\endproof

We now turn to deformations of sheaves along the twistor line.

\subsubsection{Deformation of a locally free polystable sheaf with trivial slope} Let $E_0$ be a polystable vector bundle on $S_0$ whose slope is zero, and denote by $E^\infty$ the $C^\infty-$vector bundle underlying $E_{0}$. The Kobayashi-Hitchin correspondence provides $E^\infty$ with an ASD-connection. By Lemma \ref{asd} the curvature of the connection is of $(1,1)$-type on each $S_t$. We therefore obtain holomorphic structures $E_t$ on $E^\infty$ over each $S_t$, induced by the structure $E_0$ in a canonical way. In fact we even get a holomorphic structure on $q^{*}E^\infty$; denote by $\tilde E$ the corresponding sheaf of holomorphic sections over $Z(S)$. As $E_{t}$ is holomorphic and carries an ASD-connection, it is polystable for every $t\in\mathbb{P}^{1}$. 
It is easy to see that if $E_{0}$ is stable, then $E_{t}$ is stable for every $t\in\mathbb{P}^{1}$.

\subsubsection{Deformation of an Azumaya algebra} Let now $\mathcal{A}_{0}$ be an Azumaya algebra on $S_{0}$, and let $\alpha_{0}$ be its class in $Br(S_{0})$. Choose a locally free $\alpha_{0}-$twisted sheaf $E_{0}$ such that $\mathcal{A}_{0}\simeq\mathscr{E}nd(E_{0})$. We will suppose that $E_{0}$ is $\mu_{\alpha_{0},\omega_{0}}-$stable.

The Kobayashi-Hitchin correspondence for twisted sheaves established by Wang in \cite{Wang} shows that $\mathcal{A}_{0}$ is $\mu_{\omega_{0}}-$polystable. Notice that $\mu_{\omega_{0}}(\mathcal{A}_{0})=0$, hence by section 4.4.1 the vector bundle $\mathcal{A}:=q^{*}\mathcal{A}_{0}$ carries a holomorphic structure, and for every $t\in\mathbb{P}^{1}$ its restriction $\mathcal{A}_{t}$ to the fibre $S_{t}$ is a $\mu_{\omega_{t}}-$polystable vector bundle with trivial slope. We need to show that $\mathcal{A}_{t}$ is an Azumaya algebra.

To do so, we argue as in the proof of Lemma 6.5 in \cite{Mar}, point (3). The Azumaya algebra structure on $\mathcal{A}_{0}$ is given by a holomorphic map $m_{0}:\mathcal{A}_{0}\otimes\mathcal{A}_{0}\longrightarrow\mathcal{A}_{0}$ verifying some identities among holomorphic sections. This means that $m_{0}$ is a holomorphic section of the vector bundle $\mathscr{H}om(\mathcal{A}_{0}\otimes\mathcal{A}_{0},\mathcal{A}_{0})$. But this is $\mu_{\omega_{0}}-$polystable as $\mathcal{A}_{0}$ is, hence it carries an ASD-connection, and $m_{0}$ is parallel with respect to it.

As a consequence, $m_{0}$ defines a parallel section of $\mathscr{H}om(\mathcal{A}_{t}\otimes\mathcal{A}_{t},\mathcal{A}_{t})$, hence a holomorphic map $m_{t}:\mathcal{A}_{t}\otimes\mathcal{A}_{t}\longrightarrow\mathcal{A}_{t}$. Hence $\mathcal{A}_{t}$ is an $\mathscr{O}_{S_{t}}-$algebra: as the same identities among sections which are verified on $\mathcal{A}_{0}$ are verified even on $\mathcal{A}_{t}$, it follows that $\mathcal{A}_{t}$ is an Azumaya algebra.\footnote{If $E_{0}$ is an untwisted sheaf, we can give a more direct proof. The multiplication of two holomorphic sections $\phi_1$, $\phi_2$ of $\mathcal{A}_t$ remains holomorphic (hence $\mathcal{A}_{t}$ is a sheaf of algebras on $S_{t}$): this is a consequence of the formula $\hat D(\phi_1\circ\phi_2)=\hat D\phi_1\circ\phi_2+\phi_1\circ\hat D\phi_2$, where $\hat D$ is the connection induced by $D$ on $\mathcal{A}_{0}$. 

By \cite[Thm. 1.1.6]{Cal} we just need to show that $\mathcal{A}_t$ is locally of the form $\mathscr{E}nd(F)$ for some locally free sheaf $F$ of $\mathscr{O}_{S_t}$-modules. To do so, consider the self-dual part $R_{SD}$ of the curvature $R$ of $D$. We have $R_{SD}=c \cdot Id\cdot\omega_0$ for a suitable constant $c$. By solving the equation $dd^c\phi= -\frac{c}{r}\omega_0$ on a open subset $U$, we find a holomorphic hermitian line bundle $(L,h)$ on $U$ whose curvature is $-\frac{c}{r}\omega_0$. Hence $F^{\infty}:=E_0\otimes L$ is a rank $r$ vector bundle on $U$ with a Hermite-Einstein connection, and $\mathcal{A}^\infty\cong\mathscr{E}nd^\infty(F^\infty)$ as ASD-vector bundles. Hence on $F^\infty$ we have a holomorphic structure $F_t$ compatible with the corresponding $I_t$, and $\mathcal{A}_t\cong \mathscr{E}nd(F_t)$.}

\subsubsection{Deformation of a stable twisted vector bundle} Let $\alpha_{0}\in Br(S_{0})$ and $F_{0}$ an $\alpha_{0}-$twisted locally free sheaf which is $\mu_{\alpha_{0},\omega_{0}}-$stable. Choose an $\alpha_{0}-$twisted locally free sheaf $E_{0}$ which is $\mu_{\alpha_{0},\omega_{0}}-$stable in such a way that $c_{E_{0},1}(\mathscr{F}_{0})=0$. 

We let $G_{0}:=F_{0}\otimes E_{0}^{\vee}$ and $\mathcal{A}_{0}:=E_{0}\otimes E_{0}^{\vee}$: then $\mathcal{A}_{0}$ is an Azumaya algebra, and as we saw in section 4.4.2 it is a polystable sheaf. Moreover, $G_{0}$ is a locally free sheaf of trivial slope and it has the structure of $\mathcal{A}_{0}-$module. The Kobayashi-Hitchin correspondence for twisted sheaves in \cite{Wang} shows that $G_{0}$ is a polystable sheaf.

Following section 4.4.2, $q^{*}\mathcal{A}_{0}$ is a holomorphic vector bundle, and for every $t\in\mathbb{P}^{1}$ its restriction $\mathcal{A}_{t}$ to $S_{t}$ is a polystable sheaf having the structure of Azumaya algebra. We let $\alpha_{t}$ be its class in $Br(S_{t})$.

By section 4.4.2 the polystable sheaf $G_{0}$ gives rise, for every $t\in\mathbb{P}^{1}$, to a polystable sheaf $G_{t}$ with trivial slope. The same argument used in section 4.4.2 to show that $\mathcal{A}_{t}$ is an Azumaya algebra, applied this time to $m_{t}:\mathcal{A}_{t}\otimes G_{t}\longrightarrow G_{t}$, shows that the sheaf $G_{t}$ has the structure of an $\mathcal{A}_{t}-$module.

As $G_{t}$ is an $\mathcal{A}_{t}-$module, it corresponds to an $\alpha_{t}-$twisted locally free sheaf $F_{t}$ on $S_{t}$. In particular $E_{0}$ gives rise to an $\alpha_{t}-$twisted locally free sheaf $E_{t}$ on $S_{t}$ such that $\mathscr{E}nd(E_{t})\simeq \mathcal{A}_{t}$ and $F_{t}\otimes E_{t}^{\vee}\simeq G_{t}$.

\begin{lem}
\label{lem:stabdefotwi}The sheaves $F_{t}$ and $E_{t}$ are $\mu_{\alpha_{t},\omega_{t}}-$stable.
\end{lem}

\proof We show that $E_{t}$ is $\mu_{\alpha_{t},\omega_{t}}-$stable. The proof for $F_{t}$ is similar. Suppose that $E_{t}$ is not $\mu_{\alpha_{t},\omega_{t}}-$stable, and let $\mathscr{E}_{t}\subseteq E_{t}$ in $Coh(S_{t},\alpha_{t})$ with $\mu_{E_{t},\omega_{t}}(\mathscr{E}_{t})\geq\mu_{E_{t},\omega_{t}}(E_{t})$. We suppose that $\mathscr{E}_{t}$ is $\mu_{\alpha_{t},\omega_{t}}-$stable.

We let $\mathcal{H}_{t}:=\mathscr{E}_{t}\otimes E_{t}^{\vee}$, which is an $\mathcal{A}_{t}-$module and we have $\mathcal{H}_{t}\subseteq \mathcal{A}_{t}$. The inequality $\mu_{E_{t},\omega_{t}}(\mathscr{E}_{t})\geq\mu_{E_{t},\omega_{t}}(E_{t})$ gives $\mu_{\mathcal{A}_{t},\omega_{t}}(\mathcal{H}_{t})\geq\mu_{\mathcal{A}_{t},\omega_{t}}(\mathcal{A}_{t})$, so that $\mu_{\omega_{t}}(\mathcal{H}_{t})\geq\mu_{\omega_{t}}(\mathcal{A}_{t})=0$. As $\mathcal{A}_{t}$ is $\mu_{\omega_{t}}-$polystable, this implies that $\mu_{\omega_{t}}(\mathcal{H}_{t})=0$, and that it is a direct summand of $\mathcal{A}_{t}$. In particular, it is $\mu_{\omega_{t}}-$polystable.

Using the same argument given before, the sheaf $\mathcal{H}_{t}$ gives rise to a $\mu_{\omega_{0}}-$polystable sheaf $\mathcal{H}_{0}$ on $S_{0}$, which is contained in $\mathcal{A}_{0}$, has the structure of $\mathcal{A}_{0}-$module, and $\mu_{\omega_{0}}(\mathcal{H}_{0})=\mu_{\omega_{0}}(\mathcal{A}_{0})=0$. The equivalence between $Coh(S_{0},\alpha_{0})$ and $Coh(S_{0},\mathcal{A}_{0})$ given by tensoring with $E_{0}^{\vee}$ produces then a subsheaf $\mathscr{E}_{0}$ of $E_{0}$ such that $\mu_{E_{0},\omega_{0}}(\mathscr{E}_{0})=\mu_{E_{0},\omega_{0}}(E_{0})$. But this is not possible as $E_{0}$ is $\mu_{\alpha_{0},\omega_{0}}-$stable. In conclusion, the sheaf $E_{t}$ is $\mu_{\alpha_{t},\omega_{t}}-$stable.\endproof

\subsection{Relative moduli space of twisted sheaves on twistor lines}\label{sect:twistor-space}

In this section we show that the relative moduli space of stable twisted sheaves gives us a way to deform the moduli spaces $M^{\mu}_{\alpha,w}(S,\omega)$ to an irreducible symplectic manifold (which is moreover deformation equivalent to a Hilbert scheme of points on projective K3 surface). 

We let $S$ be a K3 surface, $w=(r,0,a)\in H^{2*}(S,\mathbb{Z})$ with $r\geq 2$, $\alpha\in Br(S)$ and $\omega$ an $(\alpha,w)-$generic polarization. The K\"ahler class $\omega$ corresponds to the choice of a Riemannian metric $g$ which is compatible with the complex structure $I$ of $S$, and whose associated K\"ahler class is $\omega$. Let $\pi:Z(S)\longrightarrow\mathbb{P}^{1}$ be the twistor family of $g$: we denote $S_{t}$ the fibre of $\pi$ over $t$, which corresponds to a complex structure $I_{t}$ on $S$ associated with $t$. The metric $g$ is compatible with $I_{t}$, the associated class $\omega_{t}$ is K\"ahler, and $w$ is a Mukai vector on $S_{t}$ for every $t\in\mathbb{P}^{1}$.

Choose now a $\mu_{\alpha,\omega}-$stable $\alpha-$twisted sheaf $\mathscr{E}$ on $S$ of rank $r$, and let $E:=\mathscr{E}^{\vee\vee}$: this is a $\mu_{\alpha,\omega}-$stable $\alpha-$twisted vector bundle of rank $r$, and we let $\mathcal{A}_{0}:=\mathscr{E}nd(E)$ the corresponding Azumaya algebra. We suppose that $v_{E}(\mathscr{E})=w$. By section 4.4.2, there is holomorphic vector bundle $\mathcal{A}$ on $Z(S)$ whose restriction $\mathcal{A}_{t}$ to $S_{t}$ is an Azumaya algebra on $S_{t}$ for every $t\in\mathbb{P}^{1}$. We let $\alpha_{t}\in Br(S_{t})$ be its class and $\mathcal{A}_{t}\simeq\mathscr{E}nd(E_{t})$, where $E_t$ is the deformation of $E$ along the twistor line (see section 4.4.3).

By section 4.3.2 there is then a relative moduli space of stable twisted sheaves $p:\mathscr{M}\longrightarrow\mathbb{P}^{1}$ such that for every $t\in\mathbb{P}^{1}$ the fibre over $t$ is the moduli space $M^{\mu}_{\alpha_{t},w}(S_{t},\omega_{t})$ of $\mu_{\alpha_{t},\omega_{t}}-$stable $\alpha_{t}-$twisted sheaves whose twisted Mukai vector with respect to $E_{t}$ is $w$. 

\begin{oss}
\label{oss:univfamesiste}
On $\mathscr{M}\times_{\mathbb{P}^{1}}Z(S)$ we have a quasi-universal family: if $\mathscr{F}\in M^{\mu}_{\alpha,w}(S,\omega)$, let $F:=\mathscr{F}^{\vee\vee}$ and $\mathcal{V}_{0}:=F\otimes E^{\vee}$. We let $\mathcal{V}$ in Proposition \ref{prop:univfam} be $\mathcal{V}:=q^{*}\mathcal{V}_{0}$.
\end{oss}

We first prove some geometrical properties of the relative moduli space $p:\mathscr{M}\longrightarrow\mathbb{P}^{1}$. 

\begin{prop}
\label{prop:parteprima}Let $S$ be a K3 surface, $w=(r,0,a)\in H^{2*}(S,\mathbb{Z})$ with $r\geq 2$, $\alpha\in Br(S)$ and $\omega$ a $(\alpha,w)-$generic polarization. The relative moduli space $p:\mathscr{M}\longrightarrow\mathbb{P}^{1}$ of stable twisted sheaves verifies the following properties:
\begin{enumerate}
 \item the morphism $p$ is submersive;
 \item if $T_{p}^{*}$ denotes the relative cotangent bundle of $p$, there is a holomorphic global section of $\wedge^{2}T^{*}_{p}\otimes\mathscr{O}_{\mathbb{P}^{1}}(2)$ whose restriction to any fibre is a holomorphic symplectic form;
\end{enumerate}
\end{prop}

\proof We divide the proof in several parts.

\textit{Step 1: submersivity}. As every $\mathscr{E}\in\mathscr{M}_{t}$ is simple and the canonical bundle of a K3 surface is trivial, we have $Ext^{2}(\mathscr{E},\mathscr{E})_{0}=0$. The exact sequence (\ref{eq:extan}) implies then that the map $p$ is submersive, so that condition (1) of the statement is proved.

\textit{Step 2: section through locally free sheaves}. Let $t_{0}\in\mathbb{P}^{1}$, and choose $F\in M^{\mu}_{\alpha_{t_{0}},w}(S_{t_{0}},\omega_{t_{0}})$ a locally free sheaf. As we saw in section 4.4.3, the sheaf $F$ gives rise to a sheaf $F_{t}\in M^{\mu}_{\alpha_{t},w}(S_{t},\omega_{t})$ for every $t\in\mathbb{P}^{1}$. This produces a section 
$$s_{F}:\mathbb{P}^{1}\longrightarrow\mathscr{M},\,\,\,\,\,\,\,\,s_{F}(t):=F_{t}$$of $p$, which is holomorphic. If we let $E_{t}$ be the $\alpha_{t}-$twisted $\mu_{\alpha_{t},\omega_{t}}-$stable sheaf such that $\mathcal{A}_{t}=\mathscr{E}nd(E_{t})$ (an Azumaya algebra on $S_{t}$ whose class in $Br(S_{t})$ is $\alpha_{t}$), and $G_{t}:=F_{t}\otimes E_{t}^{\vee}$, we let $\mathcal{G}:=q^{*}G_{t}$, which is a holomorphic vector bundle on $Z(S)$. The restriction of the relative tangent bundle $T_{p}$ of $p$ to the section $s$ is $$s^{*}T_{p}\simeq R^{1}\pi_{*}\mathscr{E}nd(\mathcal{G}).$$

\textit{Step 3: relative symplectic form}. We prove that the condition (2) is fulfilled. We first notice that for every $t\in\mathbb{P}^{1}$ the restriction $T_{p|t}$ of $T_{p}$ to $\mathscr{M}_{t}$ is the tangent bundle of $M^{\mu}_{\alpha_{t},w}(S_{t},\omega_{t})$, and similarly the restriction $(T_{p}^{*})_{|t}$ of $T_{p}^{*}$ to $\mathscr{M}_{t}$ is the cotangent bundle of $M^{\mu}_{\alpha_{t},w}(S_{t},\omega_{t})$. As on $M^{\mu}_{\alpha_{t},w}(S_{t},\omega_{t})$ we have a holomorphic symplectic form (if $S_{t}$ is projective, this is done in \cite{Y3}; the proof in the general case is similar), we get an isomorphism $T_{p|t}\simeq(T_{p}^{*})_{t}$. 

This implies the existence of a line bundle $\mathscr{O}_{\mathbb{P}^{1}}(d)$ for some $d\in\mathbb{Z}$ together with an isomorphism $T_{p}\longrightarrow T_{p}^{*}\otimes p^{*}\mathscr{O}_{\mathbb{P}^{1}}(d)$. We then just need to show that $d=2$. To do so, consider a locally free sheaf $F\in\mathscr{M}_{0}$: as seen in Step 2 we have a holomorphic section $s:\mathbb{P}^{1}\longrightarrow\mathscr{M}$ of $p$, and $$s^{*}T_{p}\simeq R^{1}p_{*}\mathscr{E}nd(\mathcal{G})$$where $\mathcal{G}=q^{*}(F\otimes E_{0}^{\vee})$. By the relative Serre duality we get $$R^{1}p_{*}\mathscr{E}nd(\mathcal{G})\simeq(R^{1}p_{*}\mathscr{E}nd(\mathcal{G})^{*}\otimes K_{\pi})^{*},$$where $K_{\pi}$ is the relative canonical bundle of $\pi:Z(S)\longrightarrow\mathbb{P}^{1}$.

Now, as $\mathcal{G}$ is locally free, we have $\mathscr{E}nd(\mathcal{G})\simeq\mathscr{E}nd(\mathcal{G})^{*}$. Moreover, $K_{\pi}\simeq\mathscr{O}_{\mathbb{P}^{1}}(-2)$ (see \cite{HKLR}), hence $$R^{1}p_{*}\mathscr{E}nd(\mathcal{G})\simeq R^{1}p_{*}\mathscr{E}nd(\mathcal{G})^{*}\otimes\mathscr{O}_{\mathbb{P}^{1}}(2).$$In conclusion, $$s^{*}T_{p}\simeq s^{*}T_{p}^{*}\otimes\mathscr{O}_{\mathbb{P}^{1}}(2).$$As $s^{*}T_{p}\simeq s^{*}T^{*}_{p}\otimes\mathscr{O}_{\mathbb{P}^{1}}(d)$, it follows $d=2$. This shows that condition (2) is fulfilled.\endproof

We now prove some geometrical properties of the moduli spaces of stable twisted sheaves we are considering: in particular, we show that they are all compact and connected.

\begin{prop}
\label{prop:compconn}Let $S$ be a K3 surface, $w=(r,0,a)\in H^{2*}(S,\mathbb{Z})$ with $r\geq 2$, $\alpha\in Br(S)$ and $\omega$ a $(\alpha,w)-$generic polarization. Moreover, let $\xi$ be a representative of $\alpha$ in $H^{2}(S,\mathbb{Z})$ which is prime with $r$. The moduli space $M^{\mu}_{\alpha,w}(S,\omega)$ is a compact, connected manifold.
\end{prop}

\proof The compactness of $M^{\mu}_{\alpha,w}(S,\omega)$ is well known when $S$ is projective and a proof in the non-projective and non-twisted case has been given in \cite{T}. This proof uses in an essential way the comparison map from the moduli space of stable sheaves to the corresponding  Donaldson-Uhlenbeck compactification of the moduli space of anti-self-dual connections  in  a hermitian vector bundle on $S$. These arguments may be extended to the twisted case. We refer the reader to \cite{T} and \cite{Wang} for the ingredients.

To show that $M^{\mu}_{\alpha,w}(S,\omega)$ is connected,  we will follow the strategy used by Mukai and by Kaledin, Lehn and Sorger to prove the analogous result when $S$ is projective, $\omega$ is the first Chern class of an ample line bundle, and the sheaves are untwisted (see the proof of Theorem 4.1 in \cite{KLS}).

We first suppose that the moduli space $M^{\mu}_{\alpha,w}(S,\omega)$ is not connected, and we choose a connected component $Y$. Moreover, we fix a sheaf $F\in Y$ and a sheaf $G\in M^{\mu}_{\alpha,w}(S,\omega)\setminus Y$.

Let $p:Y\times S\longrightarrow Y$ and $q:Y\times S\longrightarrow S$ be the two projections, and consider a $p^{*}\beta\cdot q^{*}\alpha-$twisted universal family $\mathscr{F}$ on $Y\times S$. We then define two complexes $$K^{\bullet}_{F}:=\mathscr{E}xt_{p}^{\bullet}(q^{*}F,\mathscr{F}),\,\,\,\,\,\,\,\,\,\,\,K^{\bullet}_{G}:=\mathscr{E}xt_{p}^{\bullet}(q^{*}G,\mathscr{F})$$of $\beta-$twisted sheaves on $Y$. 

As the sheaves $F$ and $G$ have the same topological invariants (since their Mukai vectors are equal), letting $d:=dim(Y)$, by the Grothendieck-Riemann-Roch Theorem we have $c^{B}_{d}(K^{\bullet}_{F})=c^{B}_{d}(K^{\bullet}_{G})$, where $c^{B}_{d}$ is the component of degree $2d$ of $c^{B}$ (for some $B-$field giving the twist $\beta$).

We now compute more explicitely these twisted Chern classes, and we start from $K^{\bullet}_{G}$. We notice that if $E\in Y$, then $E$ is a stable twisted sheaf having the same slope of $G$, but which is not isomorphic to $G$. It follows that $$Ext^{0}(G,F)=Ext^{2}(G,F)=0.$$As $$\mathscr{E}xt_{p}^{j}(q^{*}G,\mathscr{F})_{E}\simeq Ext^{j}(G,E),$$it follows that $$\mathscr{E}xt_{p}^{j}(q^{*}G,\mathscr{F})=0$$if $j=0,2$, and that $\mathscr{E}xt_{p}^{1}(q^{*}G,\mathscr{F})$ is a locally free $\beta-$twisted sheaf whose rank is $d-2$.

As a consequence we have $$c^{B}_{d}(K^{\bullet}_{G})=-c^{B}_{d}(\mathscr{E}xt_{p}^{1}(q^{*}G,\mathscr{F}))=0,$$as $\mathscr{E}xt_{p}^{1}(q^{*}G,\mathscr{F})$ is a locally free $\beta-$twisted vector bundle of rank $d-2<d$: recall that $c^{B}$ of $\mathscr{E}xt_{p}^{1}(q^{*}G,\mathscr{F})$ is defined as the Chern class of some untwisted vector bundle of the same rank, hence, as this rank is smaller then the dimension of $Y$, the $d-$th $B-$twisted Chern class is trivial.

We now need to compute $c^{B}_{d}(K^{\bullet}_{F})$. To do so, we first recall that by \cite{BPS} there is locally on $Y$ a complex $$A^{\bullet}=\cdots\stackrel{a_{-1}}\longrightarrow A^{0}\stackrel{a_{0}}\longrightarrow A^{1}\stackrel{a_{1}}\longrightarrow A^{2}\longrightarrow 0$$of free sheaves such that for every $\sigma:Y'\longrightarrow Y$ and for every $j\in\mathbb{Z}$ we have $$\mathscr{E}xt^{j}_{p'}(\sigma^{*}(q')^{*}F,\sigma^{*}\mathscr{F})\simeq\mathcal{H}^{j}(\sigma^{*}A^{\bullet}),$$where $p':Y'\times S\longrightarrow Y'$ and $q':Y'\times S\longrightarrow S$ are the two projections, and where $\mathcal{H}^{j}$ denotes the $j-$th cohomology of the complex.

Let us now cover $Y$ with open subsets $U_{i}$ so that $F$ is contained in only one of them, and let us moreover suppose that the previous complex $A^{\bullet}$ exists over $U_{i}$. If $E\in U_{i}$ and $E$ is not $F$, then $\mathcal{H}^{j}(A^{\bullet})_{E}=0$, hence the rank of all the maps $a_{i}$ of the complex $A^{\bullet}$ is constant on $Y\setminus\{F\}$. But we have $$\mathcal{H}^{0}(A^{\bullet})_{F}\simeq\mathcal{H}^{0}(A^{\bullet})_{F}\simeq\mathbb{C},$$hence the rank of $a_{0}$ and $a_{1}$ at $F$ drops by 1, while the rank of $a_{i}$ is constant on $Y$ for $i\leq -1$. The same proof of Lemma 4.3 of \cite{KLS} shows that the degeneracy locus of the $a_{0}$ and $a_{1}$ is the reduced point $F$, while $a_{i}$ does not degenerate if $i\leq -1$.

Let us now consider the blow-up $\sigma:Z\longrightarrow Y$ of $Y$ at $F$ with reduced structure, and let $D$ be the exceptional divisor on $Z$. Consider the complex $$\sigma^{*}A^{\bullet}=\cdots\stackrel{\sigma^{*}a_{-1}}\longrightarrow\sigma^{*}A^{0}\stackrel{\sigma^{*}a_{0}}\longrightarrow\sigma^{*}A^{1}\stackrel{\sigma^{*}a_{1}}\longrightarrow\sigma^{*}A^{2}\longrightarrow 0.$$The degeneracy locus of $\sigma^{*}a_{0}$ and $\sigma^{*}a_{1}$ is then the exceptional divisor $D$ with reduced structure, while the $\sigma^{*}a_{i}$'s do not degenerate on $Z$ for $i\leq -1$.

The maps $\sigma^{*}a_{0}$ and $\sigma^{*}a_{1}$ hence factor through $$(A')^{0}\stackrel{a'_{0}}\longrightarrow \sigma^{*}A^{1}\stackrel{a'_{1}}\longrightarrow(A')^{2}$$where $\sigma^{*}A^{0}\subseteq (A')^{0}$, $(A')^{2}\subseteq\sigma^{*}A^{2}$, and the sheaves $$M:=(A')^{0}/\sigma^{*}A^{0},\,\,\,\,\,\,\,\,\,\,L:=\sigma^{*}A^{2}/(A')^{2}$$are supported on $D$. As in the proof of Theorem 4.1 of \cite{KLS}, Step 4, the sheaves $L$ and $M$ are characterized by canonical isomorphisms $$L\otimes\mathcal{O}_{D}\simeq Ext^{2}(F,F)\otimes\mathcal{O}_{D},\,\,\,\,\,\,\,\,\,\,\,Tor_{1}^{\mathcal{O}_{D}}(M,\mathcal{O}_{D})\simeq Hom(F,F)\otimes\mathcal{O}_{D}.$$Here the computation is done in a neighborhood of the divisor $D$.

As in \cite{KLS}, it follows from this that $$\mathscr{E}xt^{0}_{p'}(\sigma^{*}q^{*}F,\sigma^{*}\mathscr{F})\simeq\mathcal{O}_{D}(D),\,\,\,\,\,\,\,\,\,\,\mathscr{E}xt^{2}_{p'}(\sigma^{*}q^{*}F,\sigma^{*}\mathscr{F})\simeq\mathcal{O}_{D},$$viewed as $\sigma^{*}\beta-$twisted sheaves, and that $\mathscr{E}xt^{1}_{p'}(\sigma^{*}q^{*}F,\sigma^{*}\mathscr{F})$ is a locally free $\sigma^{*}\beta-$twisted sheaf of rank $d-2$. It follows that $$c^{B}_{d}(\sigma^{*}K^{\bullet}_{F})=D^{d}=-1.$$But remark that $$c^{B}_{d}(\sigma^{*}K^{\bullet}_{F})=\sigma^{*}c^{B}_{d}(K^{\bullet}_{F})=\sigma^{*}c^{B}(K^{\bullet}_{G})=0,$$getting a contradiction. In conclusion the moduli space $M^{\mu}_{\alpha,w}(S,\omega)$ has to be connected.\endproof

We can now prove the following result, which is the main result of this section, and which concludes the proof of part (1) of Theorem \ref{thm:main}:

\begin{prop}
\label{prop:diffeotwi}Let $S$ be a K3 surface, $w=(r,0,a)\in H^{2*}(S,\mathbb{Z})$ with $r\geq 2$, $\alpha\in Br(S)$ and $\omega$ a $(\alpha,w)-$generic polarization. Moreover, let $\xi$ be a representative of $\alpha$ in $H^{2}(S,\mathbb{Z})$ which is prime with $r$. Consider the relative moduli space of stable twisted sheaves $p:\mathscr{M}\longrightarrow\mathbb{P}^{1}$ along the twistor family of $(S,\omega)$.
\begin{enumerate}
 \item There is a $\overline{t}\in\mathbb{P}^{1}$ such that $\mathscr{M}_{\overline{t}}$ is an irreducible symplectic manifold which is deformation equivalent to a Hilbert scheme of points on a projective K3 surface $S$.
 \item The moduli space $M^{\mu}_{\alpha,w}(S,\omega)$ is a compact, connected complex manifold which is simply connected and carries a holomorphic symplectic form.
\end{enumerate}
\end{prop}

\proof We let $\pi:Z(S)\longrightarrow\mathbb{P}^{1}$ be the twistor family of $(S,\omega)$. By \cite[Lemma 2.1]{HSc} there is a $\overline{t}$ such that $S_{\overline{t}}$ is a projective K3 surface. The polarization $\omega_{\overline{t}}$ is $(\alpha_{\overline{t}},w)-$generic, and $w_{\xi}=v(E_{\xi})$ for some topological vector bundle $E_{\xi}$: such a topological vector bundle remains constant along $\mathbb{P}^{1}$, hence $w_{\xi}=(r,\xi,b)$ where $r$ and $\xi$ are prime to each other. It follows from Proposition \ref{prop:modtwiproj} that $M^{\mu}_{\alpha_{\overline{t}},w}(S_{\overline{t}},\omega_{\overline{t}})$ is an irreducible symplectic manifold which is deformation equivalent to a Hilbert scheme of points on $S_{\overline{t}}$.

By Proposition \ref{prop:compconn}, all the fibers are compact, connected manifolds, and by point (a) of Proposition \ref{prop:parteprima} the morphism $p$ is submersive. By the Proposition in section 1 of \cite{E}, it follows that $p$ is a smooth and proper morphism, hence it is a deformation of $M^{\mu}_{\alpha,w}(S,\omega)$, and we are done.\endproof

\subsection{Moduli spaces of locally free sheaves}

The previous results can be largely improved if we suppose something more on $M^{\mu}_{v}(S,\omega)$, namely that it parametrizes only locally free sheaves. However this case has already been considered by differential geometers. We therefore only state the following result and refer the reader to \cite{Ito} and \cite{ItNa} for the proof. 

\begin{prop}
\label{prop:locfree}Let $S$ be a K3 surface, $v=(r,\xi,a)$ a Mukai vector such that $r$ and $\xi$ are prime to each other, and $\omega$ a $v-$generic polarization. Then the open part $\mathscr{M}^{lf}$ of the relative moduli space $p:\mathscr{M}\longrightarrow\mathbb{P}^{1}$ along the twistor family of $(S,\omega)$, parameterizing locally free sheaves, is the twistor family of the moduli space $M^{\mu-lf}_{v}(S,\omega)$ of $\omega$-stable locally free sheaves with Mukai vector $v$ on $S$. 
\end{prop}

If morever $v^{2}=0$, a standard argument shows that every sheaf in $M^{\mu}_{v}(S,\omega)$ is locally free (see Remark 6.1.9 of \cite{HL}), and thus the previous proposition applies to $M^{\mu}_{v}(S,\omega)$ which is moreover compact.
The next proposition shows that compact moduli spaces of stable locally free sheaves as above may attain any even complex dimension.

\begin{prop}
\label{prop:existence}\label{Prop:dim}
Let $r$ be a positive integer, $d\in[0, 2r-2]$ be an even integer and  $g$ be an integer such that $g\le -(r^2-1)(r-1)$ and $g$ congruent to $\frac{d}{2}$ modulo $r$. 
Then there exists a K3 surface $X$ with $NS(X)$ generated by one element $\xi$ such that $\xi^2=2g-2$ and there exist torsion-free coherent sheaves $E$  on $X$ of rank $r$, $c_1(E)=\xi$ and such that $2r^2\Delta(E)-2(r^2-1)=d$. Moreover all such sheaves are locally free and irreducible. In particular they are stable with respect to any polarization on $X$ and their moduli space is a compact irreducible holomorphic symplectic manifold of dimension $d$. 
\end{prop}

\proof
 The existence of K3 surfaces $X$ with cyclic N\'eron-Severi groups was proved in \cite{LeP} whereas the existence of torsion-free sheaves $E$ with the above invariants follows from \cite[Theorem 2.7]{KuYo}. We shall check that such  sheaves are irreducible and locally free.   
Suppose
$0\to E_1\to E\to E_2\to 0$ is an exact sequence with $E_i$ coherent sheaves without torsion on $X$ of ranks $r_i$ and with $c_1(E_i)=\xi_i$, $(i=1,2)$.
 Then $\xi_1+\xi_2=\xi$, $r_1+r_2=r$ and we directly compute 
$$\Delta(E)=\frac{1}{2r}(\frac{\xi^2}{r}-\frac{\xi_1^2}{r_1}-\frac{\xi_2^2}{r_2})+\frac{r_1}{r}\Delta(E_1)+
\frac{r_2}{r}\Delta(E_2).$$
Since $g\le0$, $X$ is non-algebraic hence $\Delta(E_i)\geq0$  and thus
$$\Delta(E)\geq\frac{1}{2r}\left(\frac{\xi^2}{r}-\frac{\xi_1^2}{r_1}-\frac{\xi_2^2}{r_2}\right)=$$ $$-\frac{1}{2r_1r_2}\left(\frac{r_2\xi}{r}-\xi_2\right)^2\geq-\frac{\xi^2}{2r^2(r-1)}=\frac{1-g}{(r-1)r^2}>\frac{r^2-1}{r^2}=1-\frac{1}{r^2}.$$ 
But this implies $d>2r^2$ which contradicts our choice of $d$. Hence $E$ is irreducible. If $E$ were not locally free an easy computation would imply that the discriminant of its double dual would be negative: a contradiction to the non-algebraicity of $X$. 
\endproof

\section{The second integral cohomology}

We now study the second integral cohomology of $M^{\mu}_{v}(S,\omega)$. We will show that it carries a non-degenerate quadratic form of signature $(3,20)$, and that we have an isometry between $H^{2}(M^{\mu}_{v},\mathbb{Z})$ and $v^{\perp}$. If $M^{\mu}_{v}(S,\omega)$ is K\"ahler, it is even a Hodge isometry: as a consequence, we will show that the moduli space is projective if and only if $S$ is projective.

\subsubsection{The quadratic form}

All along this section we will let $X:=M^{\mu}_{v}(S,\omega)$ for a choice of a K3 surface $S$, a Mukai vector $v=(r,\xi,a)$ with $r$ and $\xi$ prime to each other, and a $v-$generic polarization $\omega$. We let $2n$ be its complex dimension. We start by defining a quadratic form on $H^{2}(X,\mathbb{C})$ for every holomorphic symplectic form $\sigma$ on $X$, by using the same formula as for the Beauville form of an irreducible symplectic manifold: for every $\alpha\in H^{2}(X,\mathbb{C})$, we let
 $$q_{\sigma}(\alpha):=\frac{n}{2} \int_{X}\alpha^{2}\wedge\sigma^{n-1}\wedge\overline{\sigma}^{n-1}
 \int_X\sigma^n\wedge\overline{\sigma}^n+$$ $$(1-n)\int_{X}\alpha\wedge\sigma^{n}\wedge\overline{\sigma}^{n-1}\int_{X}\alpha\wedge\sigma^{n-1}\wedge\overline{\sigma}^{n}.$$
Note that the symplectic form is always supposed to be closed so the above definition does not depend on representatives. Note also that $q_\sigma(\sigma+\overline{\sigma})=  (\int_X\sigma^n\wedge\overline{\sigma}^n)^2\neq0$ so $q_\sigma$ is non-trivial.

Recall next the definition of the "topological" quadratic form $$\tilde q_X(\alpha):=c_n\int_X\alpha^2\sqrt{\td(X)}$$where $c_n$ is a constant depending only on $n$ chosen so that the form becomes integral on $H^2(X,\Z)$ (see \cite{GHJ}, Definition 26.19 in Part III. Compact Hyperk\"ahler Manifolds). It is known that $q_\sigma$ and $\tilde q_X$ are proportional when $X$ is moreover supposed to be K\"ahler.

We finally define $\tilde H^{2,0}:=\Im((\{\tau\in H^0(\Omega^2) \ | \ \d\tau=0\}\to H^2(X, \C))$ and $\tilde h^{2,0}(X):=\dim \tilde H^{2,0}(X)$.
We first prove the following:

\begin{prop}
\label{prop:pseudobeau}Let $p:X\to C$ be a proper submersion of relative dimension $2n$ over a connected curve $C$ such that there exists a point $0\in C$ with $X_0:=p^{-1}(0)$ irreducible holomorphic symplectic. Suppose moreover that there exists a relative non-degenerate symplectic form $\sigma\in H^0(X,\Omega^2_{X/C}\otimes p^* L)$ with values in a line bundle $L$ over $C$. Let $q_t:=q_{\sigma_t}$ be the quadratic form defined by $\sigma$ on $H^2(X_t,\C)$ for each $t\in C$. Then for all $t\in C$ the quadratic form $q_{t}$ is a positive multiple of $\tilde q_{X_0}$. In particular $q_t$ is non-degenerate of signature $(3,b_{2}(X)-3)$ and 
  $\tilde h^{2,0}(X_t)=1$.
\end{prop}

\proof We may suppose that $L$ is the trivial line bundle on $C$. Indeed, for the general case it will suffice to take trivializations of $L$ over Zariski open subsets of $C$ containing $0$.

Fix some $\alpha\in H^2(X_t,\C)$ and define for $t_1,t_2\in C$:
$$q_{t_1,t_2}(\alpha):=\frac{n}{2}
\int_{X}\alpha^{2}\wedge\sigma^{n-1}_{t_1}\wedge\overline{\sigma}^{n-1}_{t_2}
\int_X\sigma^n_{t_1}\wedge\overline{\sigma}^n_{t_2}+$$ $$+(1-n)
\int_{X}\alpha\wedge\sigma^{n}_{t_1}\wedge\overline{\sigma}^{n-1}_{t_2}
\int_{X}\alpha\wedge\sigma^{n-1}_{t_1}\wedge\overline{\sigma}^{n}_{t_2}.$$
(Note again that the above formula does not depend on representatives since the syplectic forms $\sigma_{t_1}$, $\sigma_{t_2}$ are closed.)
This defines a complex function on $C\times C$ which is holomorphic in $t_1$ and antiholomorphic in $t_2$. It becomes holomorphic on $C\times C^-$, where $C^-$ denotes the curve $C$ with the opposite complex structure. Over an analytical open neighbourhood $U$ of $0$ in $C$ all fibers $X_t$ are K\"ahler. Hence for $t\in U$ the quadratic form $q_t$ is proportional to $\tilde q$. Take now $\alpha, \alpha'\in H^2(X_0,\C)$ such that $ q_0(\alpha)\neq 0$. Then the meromorphic function 
$$(t_1,t_2)\mapsto\frac{q_{t_1,t_2}(\alpha')}{q_{t_1,t_2}(\alpha)}$$
on $C\times C^-$ is constant on the diagonal $\Delta_U\subset U\times U^-\subset C\times C^-$. But $\Delta_U$ is Zariski dense in $C\times C^-$. To see this consider the system of local holomorphic curves $C_t$ on $C\times C^-$ given as images of the maps $z\mapsto (t+z,t+\bar z)$. Each curve $C_t$ passes through the reference point $(t,t)\in\Delta_U$ but its intersection with $\Delta_U$ is a piece of a "real line". Hence by the principle of isolated zeroes any holomorphic function vanishing locally on $\Delta_U$ will also vanish on the curves  $C_t$  and thus also on the three dimensional real submanifold of $C\times C^-$ they cover. Therefore the function
$(t_1,t_2)\mapsto\frac{q_{t_1,t_2}(\alpha')}{q_{t_1,t_2}(\alpha)}$ is constant on $C\times C^-$.
From this it follows that $q_t$ is proportional to $\tilde q$ for any $t\in C$.

It remains to check that $\tilde h^{2,0}(X_t)=1$ for all $t\in C$. For this we will show that the kernel $K$ of the linear map 
$$\{\tau\in H^0(X_t,\Omega^2) \ | \ \d\tau=0\}\to H^0(X_t, K_{X_t}), \ \tau\mapsto \tau\wedge\sigma^{n-1},$$
consists of $\d$-exact forms only.
Let $b_t$ be the associated bilinear form to $q_t$. Then for any $\tau\in K$ and $\alpha\in H^2(X_t,\C)$ we have
$$b_t(\tau,\alpha)=\frac{n}{2}
 \int_{X}\tau\wedge\alpha\wedge\sigma_t^{n-1}\wedge\overline{\sigma}_t^{n-1}
 \int_X\sigma_t^n\wedge\overline{\sigma}_t^n+$$ $$
 +\frac{1-n}{2}
 \int_{X}\tau\wedge\sigma^{n}_t\wedge\overline{\sigma}^{n-1}_t
 \int_{X}\alpha\wedge\sigma_t^{n-1}\wedge\overline{\sigma}^{n}_t+$$ $$
 +\frac{1-n}{2}
 \int_{X}\alpha\wedge\sigma^{n}_t\wedge\overline{\sigma}^{n-1}_t
 \int_{X}\tau\wedge\sigma_t^{n-1}\wedge\overline{\sigma}^{n}_t=0$$
 and our assertion follows since $q_t$ is non-degenerate.\endproof

\subsubsection{Isometry with $v^{\perp}$}

We now show that there is an isometry between $H^{2}(M^{\mu}_{v}(S,\omega),\mathbb{Z})$ and $v^{\perp}$ if $v^{2}>0$, and with $v^{\perp}/\mathbb{Z}\cdot v$ if $v^{2}=0$. 

We introduce some notations. If $v\in H^{2*}(S,\mathbb{Z})$, we let $v^{\perp}$ be the orthogonal of $v$ with respect to the Mukai pairing. If $v=(r,\xi,a)$ and $\xi\in NS(S)$, then the pure weight-two Hodge structure on $H^{2*}(S,\mathbb{Z})$ induces a pure weight-two Hodge structure on $v^{\perp}$: namely, a class $\alpha=(\alpha_{0},\alpha_{1},\alpha_{2})\in v^{\perp}$ is of $(1,1)-$type if and only if $\alpha_{1}\in NS(S)$.

If $\alpha=(\alpha_{0},\alpha_{1},\alpha_{2})\in H^{2*}(S,\mathbb{Q})$, we write $\alpha^{\vee}:=(\alpha_{0},-\alpha_{1},\alpha_{2})$. If $\alpha=ch(F)$ for some locally free sheaf $F$, then $\alpha^{\vee}=ch(F^{\vee})$. It is immediate to see that if $\alpha,\beta\in H^{2*}(S,\mathbb{Q})$, then $(\alpha\cdot\beta)^{\vee}=\alpha^{\vee}\cdot\beta^{\vee}$. In particular, this implies that $(\beta/\alpha)^{\vee}=\beta^{\vee}/\alpha^{\vee}$ and $(\sqrt{\alpha})^{\vee}=\sqrt{\alpha^{\vee}}$, whenever these expressions make sense.

We now introduce a morphism associating to any class in $v^{\perp}$ a rational cohomology class on the moduli space of stable (twisted) sheaves. The construction is inspired from the similar morphism which is used in the projective case (see \cite{OG}, \cite{Y1}, \cite{Mar2}, \cite{PR}). Let $\alpha\in Br(S)$, $w\in H^{2*}(S,\mathbb{Q})$ a Mukai vector and $\omega$ a $w-$generic polarization. Suppose moreover that $M^{\mu}_{\alpha,w}(S,\omega)$ is compact, and let $p:M^{\mu}_{\alpha,w}(S,\omega)\times S\longrightarrow M^{\mu}_{\alpha,w}(S,\omega)$ and $q:M^{\mu}_{\alpha,w}(S,\omega)\times S\longrightarrow S$ be the projections. 

Choosing a quasi-universal family $\mathscr{E}$ on $M^{\mu}_{\alpha,w}(S,\omega)\times S$ of similitude $\rho$ (which exists by Remark \ref{oss:univfamesiste}), we define a morphism $$\lambda_{S,\alpha,w}:w^{\perp}\longrightarrow H^{2}(M^{\mu}_{\alpha,w}(S,\omega),\mathbb{Q})$$ by letting $$\lambda_{S,\alpha,w}(\beta):=\frac{1}{\rho}[p_{*}(q^{*}(\beta^{\vee}\cdot\sqrt{td(S)})\cdot ch(\mathscr{E}))]_{1},$$where $[\cdot]_{1}$ is the part lying in $H^{2}(M^{\mu}_{\alpha,w}(S,\omega),\mathbb{Q})$. As $\beta\in w^{\perp}$, the class $\lambda_{S,\alpha,w}(\beta)$ does not depend on the chosen quasi-universal family. If $\alpha=0$ we simply write $\lambda_{S,w}$ for $\lambda_{S,0,w}$.

We now show the following, which is a generalization of known results in the projective case (see \cite{M2}, \cite{OG}, \cite{Y1}):

\begin{prop}
\label{prop:h2}Let $S$ be a K3 surface, $v=(r,\xi,a)\in H^{2*}(S,\mathbb{Z})$ where $r\geq 2$, $\xi\in NS(S)$, $(r,\xi)=1$ and $v^{2}\geq 0$. Moreover, let $\omega$ be a $v-$generic polarization. Then the image of $\lambda_{S,v}$ is contained in $H^{2}(M^{\mu}_{v}(S,\omega),\mathbb{Z})$, and
\begin{enumerate}
 \item if $v^{2}=0$, then $\lambda_{S,v}$ defines an isometry $$\overline{\lambda}_{S,v}:v^{\perp}/\mathbb{Z}\cdot v\longrightarrow H^{2}(M^{\mu}_{v}(S,\omega),\mathbb{Z});$$
 \item if $v^{2}>0$, then $\lambda_{S,v}$ is an isometry.
\end{enumerate}
\end{prop}

\proof If $v^{2}>0$, we just need to show the following properties:
\begin{enumerate}
 \item[a)] the image of $\lambda_{S,v}$ is contained in $H^{2}(M^{\mu}_{v}(S,\omega),\mathbb{Z})$;
 \item[b)] the morphism $\lambda_{S,v}$ is bijective;
 \item[c)] the morphism $\lambda_{S,v}$ is an isometry.
\end{enumerate}
Let $\mathscr{E}$ be a quasi-universal family of similitude $\rho$ on $M^{\mu}_{v}(S,\omega)\times S$, and fix a locally free $\mu_{\omega}-$stable vector bundle $F$ of Mukai vector $v$. Let $w:=v_{F}(F)=(r,0,a-\xi^{2}/2r)$ and $$f:M^{\mu}_{v}(S,\omega)\longrightarrow M^{\mu}_{0,w}(S,\omega),\,\,\,\,\,\,\,\,f(\mathscr{F}):=\mathscr{F}\otimes F^{\vee}$$which is an isomorphism (see Remark \ref{oss:isotwist}). 

We let $q:M_{0,w}^{\mu}(S,\omega)\times S\longrightarrow S$ be the projection, and $$\mathscr{E}':=(f\times id_{S})_{*}\mathscr{E}\otimes q^{*}F^{\vee},$$which is a quasi-universal family of similitude $\rho$ on $M^{\mu}_{0,w}(S,\omega)\times S$. Moreover, as $f$ is an isomorphism, the morphism $$f_{*}:H^{2}(M^{\mu}_{v}(S,\omega),\mathbb{Z})\longrightarrow H^{2}(M^{\mu}_{0,w}(S,\omega),\mathbb{Z})$$is easily checked to be an isometry.

Now, we let $$h:H^{2*}(S,\mathbb{Z})\longrightarrow H^{2*}(S,\mathbb{Q}),\,\,\,\,\,\,\,h(\beta):=\frac{\beta\cdot ch(F^{\vee})}{\sqrt{ch(F\otimes F^{\vee})}}.$$We let $(\cdot,\cdot)_{S}$ be the Mukai pairing on $S$ and $[\cdot]_{2}$ the part lying in $H^{4}(S,\mathbb{Q})$. If $\beta\in v^{\perp}$ we have $$(h(\beta),w)_{S}=-\bigg[\frac{\beta^{\vee}\cdot ch(F)}{\sqrt{ch(F\otimes F^{\vee})}}\cdot v_{F}(F)\bigg]_{2}=$$ $$=-[\beta^{\vee}\cdot ch(F)\cdot\sqrt{td(S)}]_{2}=(\beta,v)_{S}=0,$$so that $$h:v^{\perp}\longrightarrow w^{\perp}.$$The same argument shows that it is an isometry. We even have $f_{*}(\lambda_{S,v}(\beta))=\lambda_{S,w}(h(\beta))$. Indeed $$f_{*}(\lambda_{S,v}(\beta))=\frac{1}{\rho}[f_{*}p_{*}(q^{*}(\beta^{\vee}\sqrt{td(S)})\ch(\mathscr{E}))]_{1}=$$ $$=\frac{1}{\rho}[p_{*}((f\times id_{S})_{*}q^{*}(\beta^{\vee}\sqrt{td(S)})ch(\mathscr{E}'))]_{1}=$$ $$=\frac{1}{\rho}[p_{*}(q^{*}(h(\beta)^{\vee}\sqrt{td(S)})ch(\mathscr{E}'))]_{1}=\lambda_{S,w}(h(\beta)).$$In conclusion, we see that $\lambda_{S,v}$ verifies the properties a), b) and c) above if and only if $\lambda_{S,w}$ verifies them.

Now, consider the twistor line of $(S,\omega)$ and let $p:\mathscr{M}\longrightarrow\mathbb{P}^{1}$ be the associated relative moduli space. As we can define $\lambda_{S,v}$ in a relative way using relative quasi-universal families (which exist by Remark \ref{oss:univfamesiste}), properties a), b) and c) above are verified on a fibre if and only if they are verified all along the twistor line. It follows that $\lambda_{S,w}$ verifies a), b) and c) if and only $\lambda_{S_{t},w_{t}}$ verifies them for some $t\in\mathbb{P}^{1}$. 

As we saw before, there is $t$ such that $S_{t}$ is projective, and in this case $\lambda_{S_{t},w_{t}}$ is an isometry by \cite{Y3}, hence we are done. If $v^{2}=0$, the proof is similar: the only difference is about the fact that $\mathbb{Z}\cdot v$ is the kernel of $\lambda_{S,v}$, which holds in the general case as it holds over a projective K3 surface (see \cite{M2}).\endproof

An immediate corollary of the previous Proposition is the following:

\begin{cor}
\label{cor:hodge}Let $S$ be a K3 surface, $v=(r,\xi,a)\in H^{2*}(S,\mathbb{Z})$ where $\xi\in NS(S)$, $r\geq 2$, $(r,\xi)=1$ and $v^{2}\geq 0$. If $\omega$ is a $v-$generic polarization and $M^{\mu}_{v}(S,\omega)$ is K\"ahler, then the morphism $\lambda_{v}$ is a Hodge isometry.
\end{cor}

Theorem \ref{thm:proj} can now be seen as a corollary of the previous results:

\begin{cor}
\label{cor:proj}Let $S$ be a K3 surface, $v=(r,\xi,a)\in H^{2*}(S,\mathbb{Z})$ where $\xi\in NS(S)$, $r\geq 2$, $(r,\xi)=1$ and $v^{2}\geq 0$. If $\omega$ is a $v-$generic polarization, then $M^{\mu}_{v}(S,\omega)$ is projective if and only if $S$ is projective.
\end{cor}

\proof First, notice that if $S$ is projective, then $M^{\mu}_{v}(S,\omega)$ is projective by Theorem \ref{thm:change}. 

Suppose now that $S$ is not projective, we want to prove that $M^{\mu}_{v}(S,\omega)$ is not projective as well. Suppose that $M^{\mu}_{v}(S,\omega)$ is projective: in particular this implies that it is K\"ahler, hence by part (1) of Theorem \ref{thm:main} it follows that it is an irreducible symplectic manifold. 

Recall that an irreducible symplectic manifold $X$ is projective if and only if there is a line bundle $L$ on $X$ such that $q(L)>0$, where $q$ is the Beauville form of $X$ (see \cite{Hu99}). Hence there is a line bundle $L$ on $M^{\mu}_{v}(S,\omega)$ such that $q(L)>0$, where $q$ is the Beauville form on $M^{\mu}_{v}(S,\omega)$, which coincides with the non-degenerate quadratic form we defined in the previous section.
 
Moreover, by Corollary \ref{cor:hodge}, as $M^{\mu}_{v}(S,\omega)$ is K\"ahler we have that $\lambda_{v}$ is a Hodge isometry. There is then $\alpha\in v^{\perp}$ of type $(1,1)$ (with respect to the Hodge structure on $v^{\perp}$) such that $\lambda_{v}(\alpha)=c_{1}(L)$, and $(\alpha,\alpha)_{S}>0$.

Let us now describe $v^{\perp}\otimes\mathbb{Q}$. First, an element $(0,\zeta,b)\in\widetilde{H}(S,\mathbb{Q})$ is in $v^{\perp}\otimes\mathbb{Q}$ if and only if $b=\zeta\cdot\xi$. As $(0,\zeta,\zeta\cdot\xi)=e^{\xi/r}\cdot(0,\zeta,0)$, we have $$e^{\xi/r}\cdot H^{2}(S,\mathbb{Q})\subseteq v^{\perp}.$$

It is easy to see that $e^{\xi/r}\cdot(2r^{2},0,v^{2})\in v^{\perp}\otimes\mathbb{Q}$, hence $$e^{\xi/r}\cdot\mathbb{Q}(2r^{2},0,v^{2})\subseteq v^{\perp}\otimes\mathbb{Q}.$$This implies that $$v^{\perp}\otimes\mathbb{Q}=e^{\xi/r}\cdot(H^{2}(S,\mathbb{Q})\oplus\mathbb{Q}(2r^{2},0,v^{2})),$$so that the $(1,1)-$part $(v^{\perp})^{1,1}$ of $v^{\perp}\otimes\mathbb{Q}$ is $$(v^{\perp})^{1,1}=e^{\xi/r}\cdot(NS_{\mathbb{Q}}(S)\oplus\mathbb{Q}(2r^{2},0,v^{2})),$$where $NS_{\mathbb{Q}}(S):=NS(S)\otimes\mathbb{Q}$.

The direct sum is orthogonal with respect to the Mukai pairing, and it is easy to see that $$(e^{\xi/r}(2r^{2},0,v^{2}))^{2}=-4r^{2}v^{2}\leq 0,$$as $v^{2}\geq 0$. Moreover, as $S$ is non-projective the lattice $e^{\xi/r}NS_{\mathbb{Q}}(S)$ is negative semi-definite. It follows that $(v^{\perp})^{1,1}$ is negative semi-definite, hence for every $\alpha\in(v^{\perp})^{1,1}$ we have $(\alpha,\alpha)_{S}\leq 0$, which is not possible. In conclusion, if $S$ is not projective, the moduli space cannot be projective, and we are done.\endproof

\par\bigskip
\par\bigskip

Institut \'Elie Cartan, UMR 7502, Universit\'e de Lorraine, CNRS, INRIA, Boulevard des Aiguillettes, B.P. 70239, 54506 Vandoeuvre-l\`es-Nancy Cedex, France 

\email{Arvid.Perego@univ-lorraine.fr}

\email{Matei.Toma@univ-lorraine.fr}
\end{document}